\documentclass[11pt]{article}
\usepackage{graphicx}
\usepackage{amsmath,amsthm}
\usepackage{amssymb}
\usepackage[latin1]{inputenc}
\usepackage{enumerate}
\usepackage{slashed}
\usepackage{pb-diagram}
\usepackage{slashed}
\usepackage{color}
\theoremstyle{plain}
\newtheorem{theo}{Theorem}

\newtheorem{obs}{Observation}
\newtheorem{assumption}{Assumption}
\newtheorem{define}[theo]{Definition}
\newtheorem{lemma}[theo]{Lemma}

\newcommand{\End}{\mathrm{End}}
\newcommand{\ord}{\mathrm{ord}_G}

\newcommand{\tr}{\mathrm{tr}}

\newcommand{\diver}{\mathrm{div}}
\newcommand{\str}{\mathrm{Str}}
\newcommand{\vol}{\mathrm{vol}}
\newcommand{\Tr}{\mathrm{Tr}}
\newcommand{\rmx}{\mathrm{x}}
\newcommand{\cE}{\mathcal{E}}

\newcommand{\cD}{\mathcal{D}}

\newcommand{\frR}{\mathfrak{R}}

\newcommand{\rke}{\mathrm{rk}(E)}
\newcommand{\dotg}{\dot{g}}

\newcommand{\grad}{\mathrm{grad}}

\newcommand{\cA}{\mathcal{A}}
\newcommand{\cR}{\mathcal{R}}
\newcommand{\cow}{\overline{w}}

\newcommand{\cO}{\mathcal{O}}
\newcommand{\bbC}{\mathbb{C}}
\newcommand{\cL}{\mathcal{L}}

\newcommand{\eps}{\epsilon}

\newcommand{\frP}{\mathfrak{P}}
\newcommand{\caL}{\mathcal{L}}
\newcommand{\tie}{\tilde{e}}
\newcommand{\oi}{\overline{i}}
\newcommand{\oj}{\overline{j}}
\newcommand{\bw}{\bar{w}}

\newcommand{\frB}{\mathfrak{B}}
\newcommand{\frD}{\mathfrak{D}}

\newcommand{\bbR}{\mathbb{R}}

\newcommand{\bbZ}{\mathbb{Z}}

\newcommand{\ome}{\omega}
\newcommand{\bome}{\bar{\omega}}

\newcommand{\p}{\partial}
\renewcommand{\Re}{\mathrm Re}
\newcommand{\strr}{\mathrm{str}}
\newcommand{\Str}{\mathrm{Str}}
\newcommand{\oz}{\overline{z}}
\title{\textbf{Asymptotics of a Modified Holomorphic Analytic Torsion}}
\author{Andr\'es Larra\'in-Hubach\footnote{Department of Mathematics, University of Dayton, 300 College Park Dayton, OH 45469.  alarrainhubach1@udayton.edu.}}
\date{}
\begin{document}
\maketitle

\begin{abstract}
We prove a formula for the first few terms of the asymptotic expansion of the holomorphic analytic torsion of the Dirac operator modified by the Clifford action of a real three-form.
\end{abstract}

\leftline{{\bf MSC 2020:} 58J20, 58J99.}
\tableofcontents
\section{Introduction}

The purpose of this paper is to analyze the asymptotic behavior of a modification of the holomorphic analytic torsion of Ray-Singer and obtain  results similar to those of \cite{bis2} and \cite{fin}. 

Let $(M,J,g^{TM})$ be a K\"ahler manifold with complex structure $J$, $L$ a \textbf{positive} holomorphic line bundle over $M$, and $E$ a holomorphic bundle over $M$. Both $L$ and $E$ are equipped with their Chern connections denoted by $\nabla^L$ and $\nabla^E$ respectively. Using the canonical Spin$^c$ structure induced by the complex structure, one defines the Dirac operator $D_p$ acting on $\Gamma(\Lambda^{0,*}T^*M\otimes L^p\otimes E)$, where $L^p$ denotes the $p$-tensor power of $L.$ In this setting, Bismut and Vasserot \cite{bis2} proved an explicit formula for the leading term of the holomorphic analytic torsion of $D_p$ as $p\to\infty$. The subleading terms of this expansion were computed, under additional hypotheses on the line bundle $L$, by Finski in \cite{fin}.

Here, we consider a similar situation but change $D_p$ to $\frD_p=D_p+c(B)$, where $c(B)$ denotes Clifford multiplication by a real closed $3$-form $B$. The new operator is still odd with respect to the canonical $\bbZ_2$-grading induced by the spin$^c$ structure. This type of modification of $D_p$ has already been used in the context of non-commutative cocycles in \cite{reza} and index theory in \cite{bis}.

Mathai and Wu  defined a similar twisting for the analytic torsion \cite{mw1}, \cite{mw2}. They used an operator $D^E=\nabla^E+H\wedge\cdot$, where $\nabla^E$ is a flat connection on a bundle $E$ and $H$ is an odd-degree closed differential form. The operator $D^E$ acts on $E$-valued differential forms. Our main motivation was to obtain analogous results  for a special  twisting in the spin$^c$ case, but using the holomorphic analytic torsion instead of the analytic one. For twistings of the Capell-Miller holomorphic torsion see \cite{huang}.
 
Denoting by $\theta_p(z)$ the $\theta$-function associated to $\frD_p$ (see definition 2 below), our main result is the following
\begin{theo}
As $p\to\infty$, there is an asymptotic expansion
\begin{equation}\label{introthm}
\theta_p'(0)=\sum_{i\geq 0}(\alpha_i\ln p+\beta_i)p^{n-i}
\end{equation}
such that
\begin{equation}\label{main}
\begin{split}
\alpha_0&=\frac{n\,\rke}{2}\int_M 
\mathrm{det}(\dot{R}^L/2\pi)\,dv_y\\
\beta_0&=\frac{\rke}{2}\int_M\mathrm{det}(\dot{R}^L/2\pi) \ln(\mathrm{det}(\dot{R}^L/2\pi))\,dv_y\\
\end{split}
\end{equation}
Under the additional hypothesis for the curvature of $L$ that
\begin{equation}\label{constraint}
\ome=\frac{\sqrt{-1}}{2\pi}R^L(\cdot,\cdot)=g^{TM}(J\cdot,\cdot)
\end{equation}
we also compute
\begin{equation}\label{main1}
\begin{split}
\alpha_1&=\frac{n \rke}{2\pi} \int_M|\frB|^2_{y}dv+\frac{n }{2}\int_Mc_1(E)\frac{\ome^{n-1}}{(n-1)!}\\
&\,\,\,+ \frac{(3n+1)\rke}{12} \int_Mc_1(TM)\frac{\ome^{n-1}}{(n-1)!}\\
&\,\,\,-\int_M\big[\dfrac{\mathrm{str} (N\mathfrak{H}_{j\oj}e^{-2\pi u N})}{(1-e^{-2\pi u})^{n+1}}({u}+{u}e^{-2\pi u}-\frac{1}{\pi}(1-e^{-2\pi u}))\big)\big]_{u=0}\,dv,
\end{split}
\end{equation}
where $\mathfrak{H}$ is a certain two-form defined in \eqref{curvrel5} that depends on  $\frB$,  and $\int_M[\cdot]_{u=0}\,dv$ means taking the $u=0$ coefficient.
\end{theo}

Notice that  the leading terms of the asymptotic expansion are the same as in the untwisted case \cite{bis2}. The proof that we present of \eqref{main} in the upcoming sections is based on the modified Getzler's grading technique explained in \cite[sect 3.]{alh}. The proof of \eqref{main1} is adapted from the argument in \cite[Sec. 4]{fin} and it is based on Getzler's rescaling technique. It is interesting to see the grading and  rescaling techniques, both ubiquitous in index theory, used in the same calculation.


For future research, we would like to understand what an anomaly formula would look like in the case of the torsion of $\frD_p$. This seems to be a difficult question because the familiar rescaling techniques used for example in \cite{ly} do not work in this situation. 

\textbf{Acknowledgements.} The author would like to thank Professor Siarhei Finski who pointed out that the theory developed in \cite{fin} could be easily modified to obtain the formula for \eqref{main1} above.
\section{Holomorphic Torsion modified by a Three-Form}
Let $M$ be a closed and complex manifold of complex dimension $n$ with complex structure $J$ and K\"ahler metric $g^{TM}$. The corresponding K\"ahler form is denoted by
\begin{equation}\label{kahlerform}
\tilde{\Theta}=g^{TM}(J\cdot,\cdot)=\langle J\cdot,\cdot\rangle
\end{equation}
It's important to mention that we use both notations $g^{TM}$ and $\langle \rangle$ for the metric and its $\bbC$-bilinear extension throughout the paper. The volume form is given by $dv=\frac{\tilde{\Theta}^n}{n!}.$

The metric induces a canonical spin$^c$-structure on $M$. The corresponding  spin bundles over $M$ can be identified as $S^+=\Lambda^{0,\mathrm{Even}}T^*M$ and $S^-=\Lambda^{0,\mathrm{Odd}}T^*M$. The Levi-Civita connection $\nabla^{TM}$ induces a connection on the spin bundle that coincides with the connection $\nabla^{\Lambda^{0,*}}$ on the bundle $\Lambda^{0,*}T^*M$. 

The underlining real manifold has real dimension $2n.$ We denote by $\{w_j\}_{j=1}^n$ an orthonormal frame of $T^{1,0}M$, while $\{\bw_j\}_{j=1}^n$ denotes the corresponding orthonormal frame of $T^{ 0,1}M$. Using these frames, we define a real orthonormal frame $\{e_k\}_{k=1}^{2n}$ for $TM$  by
\begin{equation}\label{realframe}
\begin{split}
e_{2j-1}&=\frac{1}{\sqrt{2}}(w_j+\bw_j)\\
e_{2j}&=\frac{\sqrt{-1}}{\sqrt{2}}(w_j-\bw_j)\\
\end{split}
\end{equation}
We use Einstein's summation convention. That is, repeated indices are summed over.

\subsection{The Modified Dirac Operator}
Given a real vector $v\in TM$, its Clifford action on $S^+\oplus S^-$ is defined by
$$c(v)=\sqrt{2}\big(\bar{v}^{(1,0),\star}\wedge -\iota_{v^{(0,1)}}\big),$$
where $v=v^{(1,0)}+v^{(0,1)}\in T^{(1,0)}M\oplus T^{(0,1)}M$ and $\bar{v}^{(1,0),\star}$ is the metric dual of $v^{(1,0)}.$

We  use two holomorphic bundles over $M$. The first one is denoted by $E$  with hermitian metric $h^E$ and Chern connection $\nabla^E.$ The second one is a line bundle $L$ with metric $h^L$ and Chern connection $\nabla^L$.  The curvature $R^L$ induces an endomorphism $\dot{R}^L$ of $T^{(1,0)}M$ defined by $R^L(X,\bar{Y})=g(\dot{R}^LX,\bar{Y})$ for $X,Y\in T^{(0,1)}M.$ We denote 
\begin{equation}\label{omegaddef}
\omega_d=-\sum_{l,m}R^L(w_l,\bar{w}_m)\bar{w}^m\wedge \iota_{\bar{w}_l}. 
\end{equation}
Given a natural number $p\geq 1$, the $p$-tensor power of $L$ is denoted by $L^p.$
For the arguments in this paper, we also assume that $L$ is a positive line bundle. In other words, the curvature $R^L$ associated to $h^L$ is such that $\sqrt{-1}R^L$ is a positive $(1,1)$-form. We use the notation $\omega=\frac{\sqrt{-1}}{2\pi}R^L.$

Let $\cE_p=\Lambda^{0,*}\otimes E\otimes L^p=(S^+\oplus S^-)\otimes E\otimes L^p$, the connections on each bundle induce a connection on $\cE$ which we denote by $\hat{\nabla}^{\cE_p}.$ It follows that $\cE_p$ is $\bbZ_2$-graded with $\cE_p^\pm=S^\pm\otimes E\otimes L^p.$

The corresponding Dirac operator that acts on smooth sections of $\cE_p$ is defined by
\begin{equation}\label{dirac}
D^{\cE_p}=\sum_{j=1}^{2n}c(e_j)\hat{\nabla}^{\cE_p}_{e_j}.
\end{equation}
This is an odd operator with respect to the $\bbZ_2$-grading defined above.

Let $B\in \Omega^3(M)$ be a \textbf{closed} real three-form such that $c(B)^\star=c(B)$. The twisted Dirac operator  is defined by
\begin{equation}\label{tdirac}
\frD_p=D^{\cE_p}+c(B)
\end{equation}
This is still a self-adjoint and odd operator. 
\begin{define}
We  use  $B$ and $\nabla^{\Lambda^{0,*}}$ to define a new connection on the spin bundle given by
\begin{equation}\label{bnabla}
\nabla^S_X=\nabla^{\Lambda^{0,*}}+c(\iota_XB),
\end{equation}
for any vector $X\in TM.$
\end{define}
The new connection $\nabla^S$ induces a new connection on $\cE_p$ denoted by $\nabla^{\cE_p}$.  Using this new connection, we get a familiar formula for the square of $\frD_p$ given by
\begin{equation}\label{tdirac2}
\frD_p^2=(\nabla^{\cE_p})^\star \nabla^{\cE_p}+\frac{r^M}{4}+c(R^E+\frac{1}{2}R^{\det}+pR^L)-2|B|^2,
\end{equation}
where $r^M$ is the scalar curvature of $M$, $R^{\det}$ is the curvature of the induced connection on the determinant bundle of $S$ by the Levi-Civita connection, and $|B|$ is the pointwise norm of the closed three-form $B$. 

The following spectral gap is crucial to this paper.
\begin{lemma}
Given $p$ large enough, there is a  constant $C>0$ such that 
\begin{equation}\label{gap}
\mathrm{Spec}(\frD_p^2)\subset \{0\}\cup (Cp,\infty).
\end{equation}
Also, $\mathrm{Ker}(\frD^-_p)=\{0\}$ and $\mathrm{Ker}(\frD_p^+)\subset \Gamma(\Lambda^{(0,0)}\otimes E\otimes L^p).$
\end{lemma}
\begin{proof}
 This  is a consequence of \eqref{tdirac2}, the positivity of $L$, and the compactness of $M$. See \cite[thm. 1.5.8.]{mm}. 
 \end{proof}
 \subsection{The Modified Torsion}
 Given a vector bundle endomorphism $A:\cE_p\to\cE_p$, $\strr(A_y)$ denotes its supertrace over the fiber $(\cE_p)_y$ for $y\in M$, and $\str(A)=\int_M \strr(A_y)\,dv_y$. Here, we integrate against the volume form $dv_y$.

Let $\Pi_p$ denote the $L^2$-orthogonal projection onto $\mathrm{Ker}\, \frD_p$ and $\Pi_p^\perp=1-\Pi_p.$ The $\theta$-function of $\frD_p$ is defined by
\begin{equation}\label{zeta}
\theta_p(z)=-\frac{1}{\Gamma(z)}\int_0^\infty u^{z-1}\Str(Ne^{-\frac{u}{2}\frD_p^2}\Pi_p^\perp) \,du,
\end{equation}
where $N$ is the number operator which acts on $\Lambda^{0,k}T^*M$ by multiplication by $k$. It is well known that $\theta_p(z)$ is a meromorphic function which is holomorphic at $z=0$. 
\begin{define}
The twisted analytic torsion  is defined  by $e^{-\frac{1}{2}\theta_p'(0)}$. See \cite[def. 5.5.4]{mm}
\end{define}

For the proof of \eqref{main}, two useful observations are needed: first, the spectral gap \eqref{gap} implies that $\Str(N\Pi_p)=0$ for $p$ large enough and, second, a simple substitution modifies the definition of $\theta_p(z)$ to

\begin{equation}\label{theta}
p^z\theta_p(z)=-\frac{1}{\Gamma(z)}\int_0^\infty u^{z-1}\Str(Ne^{-\frac{u}{2p}\frD_p^2}\Pi_p^\perp) \,du,
\end{equation}

\section{Expansions} 
In this section we obtain two asymptotic expansions of the heat kernel $e^{-\frac{u}{2p}\frD_p}$. One in powers of $p$, valid for $p$ large, and the other one in terms of powers of $u$ which is valid as $u\to 0^+.$
\subsection{Heat Kernel}
Let $y\in M$ be the point that we use as the center of our fixed normal coordinate chart. We use parallel transport, with respect to the connection  $\nabla^S$, along geodesics emanating from $y$ to extend a frame of $S_y$ to the whole coordinate chart. Notice that we do not use the Levi-Civita connection to build this local frame. Similarly, using parallel transport, we extend  orthonormal frames of $E_y$ and $L_y$ to synchronous frames over the whole coordinate chart. From this point on, $\{\rmx_j\}_{j=1}^{2n}$  denote the normal local coordinates where $y$ corresponds to $\rmx=0.$

In these local coordinates (see \eqref{gammaexpgen}) we have that
\begin{equation}\label{nablaloc}
\nabla^{\cE_p}_{\p_{\rmx_i}}=\p_{\rmx_i}+\beta_{ij}\rmx_j+p\Gamma_i^L+\Gamma_i^E,
\end{equation}
where $\beta_{ij}\in\mathrm{End}(S_y).$ The connection symbols $\Gamma^L$ and $\Gamma^E$ vanish at $\rmx=0.$

Define
\begin{equation}\label{flathk}
q_t(x,y)=\frac{1}{(4\pi t)^n}e^{\frac{-d^2(x,y)}{4t}},
\end{equation}
where $d(x,y)$ is the geodesic distance between $x$ and $y$ in $M.$ We denote the integral kernel of $e^{-t\frD_p^2}$ by $e^{-t\frD_p^2}(x,y).$

The classical theory developed in \cite[ch. 7]{roe} or \cite[ch. 3]{ros} implies that there is an asymptotic expansion
\begin{equation}\label{exp1}
e^{-t\frD_p^2}(x,y)\sim_{t\to 0^+} \chi(d(x,y))q_t(x,y)\sum_{j=0}^\infty t^j\Theta_{j,p}(x,y).
\end{equation}
Here $\chi$ is a smooth cut-off equal to $1$ for $d(x,y)\leq \frac{\mathrm{inj^2}}{4}$, where $\mathrm{inj}$ denotes the injectivity radius of $M$. This expansion is valid in the $C^l$-norm of $M\times M$ for any $l>0.$ It is also uniform in $t$ if $t$ takes values inside a compact subset of $(0,\infty).$ See \cite[sec. 2]{alh}.

Using the coordinates $\{\rmx_j\}$, we can write $\Theta_{j,p}(\rmx)$ instead of $\Theta_{j,p}(x,y)$. The $\Theta_{j,p}(\rmx)$ are smooth functions supported near $\rmx=0$ and defined by the following recurrence relation
\begin{equation}\label{rec}
\begin{split}
\Theta_{0,p}(\rmx)&=C|g|^{-1/4}(\rmx)\\
\Theta_{j,p}(\rmx)&=\frac{-1}{|g|^{1/4}}\int_0^1s^{j-1}|g|^{1/4}(s\rmx)(\frD^2_p\Theta_{j-1,p})(s\rmx)\,ds,
\end{split}
\end{equation}
where $|g|$ denotes the determinant of the metric $g$ and we integrate along geodesics emanating from $y$, which corresponds to $\rmx=0$.

From \eqref{exp1} we get the following expansion on the diagonal
\begin{equation}\label{exp2}
p^{-n}e^{-\frac{u}{2p}\frD_p^2}(y,y)\sim 
\frac{1}{(2\pi)^n}\sum_{j=0}^\infty \frac{u^{j-n}}{2^jp^j}\Theta_{j,p}(y,y)=\frac{1}{(4\pi p)^n}\sum_{j=-n}^\infty \Theta_{j+n,p}(y)\frac{u^j}{2^jp^j}
\end{equation}

Since convergence in \eqref{exp2} is uniform in $y\in M$, we also get
\begin{equation}\label{expint}
p^{-n}\Str(Ne^{-\frac{u}{2p}\frD_p^2})\sim \frac{1}{(2\pi)^n}\sum_{j=0}^\infty \frac{u^{j-n}}{2^jp^j}\int_M \strr(N\Theta_{j,p}(y,y))\,dv_y
\end{equation}
Using the local coordinates defined above, each $\Theta_{j,p}(\rmx)$ can be expanded as a formal sum of the form $\sum_{I,l}a_{I,l} \rmx_Ip^l$ with the $a_{I,l}\in \mathrm{End}(S_y\otimes E_y\otimes L^p_y)$ independent of $p$ and $\rmx$. In these expansions, for a multi-index $I=(i_1,\ldots,i_k)$, we define $\rmx_I=\rmx_{i_1}\ldots\rmx_{i_k}.$ We  need to know the largest possible power of $p$  in each $\Theta_{j,p}(y,y)=\Theta_{j,p}(0).$

\begin{define}\label{pgrad}
Consider the vector space of differential operators written formally as sums 
\begin{equation}\label{generalpde}
T=\sum_{l,I,J}a_{l,I,J}p^l\rmx_J\p_{\rmx}^J,
\end{equation}
where the $a_{l,I,J}\in \mathrm{End}(S_y\otimes E_y\otimes L^p_y)$ are independent of $p$. The $p$-grading assigns the following values: $\mathrm{deg}(\rmx_j)=-1$, $\mathrm{deg}(p)=2$ and $\mathrm{deg}(\p_{\rmx_j})=1.$ The coefficients $a_{l,I,J}$ have degree zero. The $p$-grading of a sum is the maximum of the $p$-gradings of its summands. 
\end{define}
The most important property of the $p$-grading is the following
\begin{lemma}
 Given two differential operators $T$ and $S$ as in \eqref{generalpde}, one has $\mathrm{deg}(T\circ S)\leq \mathrm{deg}(T)+\mathrm{deg}(S)$.
 \end{lemma}
 \begin{proof}
 See  \cite[thm. 3]{alh}.
\end{proof}
It follows from \eqref{nablaloc} and definition \ref{pgrad} above that $\mathrm{deg}(\nabla^{\cE_p}_{\p_{\rmx_i}})\leq 1$, $\mathrm{deg}(\frD_p)\leq 1,$ and  $\mathrm{deg}(\frD_p^2)\leq 2$. The recurrence \eqref{rec} implies that $\mathrm{deg}(\Theta_{0,p})\leq 0$ and $\mathrm{deg}(\Theta_{j,p})\leq 2j.$ Notice that the highest possible power of $p$ that appears in $\Theta_{j,p}(0)$ is $p^j$. This follows from \eqref{rec} and the observation that a summand in $\Theta_{j,p}$ with a power of $p$ higher than $j$, that does not vanish at $\rmx=0$, would have $\mathrm{deg}>2j.$ 

This implies that we have an expansion
\begin{equation}\label{pexp1}
\Theta_{l,p}(y)=\sum_{i\geq 0}\theta_{l,p,i}(y)p^{l-i},
\end{equation}
where the $\theta_{l,p,i}(y)$ are independent of $p.$
Replacing in \eqref{exp2} we get
\begin{equation}\label{pexp2}
e^{-\frac{u}{2p}\frD_p^2}(y,y)\sim \frac{1}{(4\pi)^n}\sum_{j\geq -n}\sum_{i\geq 0}\theta_{j+n,p,i}(y)\frac{p^{n-i}u^j}{2^j}
\end{equation}
This gives, comparing with \cite[Thm. 2.5]{fin}, that
\begin{equation}\label{pexp3}
e^{-\frac{u}{2p}\frD_p^2}(y,y)\sim \sum_{i\geq 0}a_{i,u}(y)p^{n-i}=\sum_{i\geq 0}\big(\frac{1}{(4\pi)^n}\sum_{j\geq -n}\theta_{j+n,p,i}(y)\frac{u^j}{2^j}\big)p^{n-i}
\end{equation}
Similarly, \cite{fin} uses an expansion, valid as $u\to 0^+$, of each $a_{i,u}(y)$ given by
\begin{equation}\label{auxexp}
\begin{split}
a_{i,u}(y)&=\sum_{j\geq -n}a_i^{[j]}(y)u^j\\
&=\sum_{j\geq -n}\big( \frac{1}{(4\pi)^n}2^{-j}\theta_{j+n,p,i}(y)\big)u^j
\end{split}
\end{equation}
\subsection{Expansion of Analytic Torsion}
We use \eqref{pexp3} and \eqref{auxexp} to expand $\theta_p'(0)$ in powers of $p>>0$. We need the following well-known fact
$$\frac{1}{\Gamma(z)}=z-\Gamma'(1)z^2+\cO(z^3),\,\,\,\,\,\mathrm{as}\,\,z\to 0.$$
Start differentiating both sides of \eqref{theta}
\begin{equation}\label{thetapexp}
\begin{split}
\theta_p(0)\ln p+\theta_p'(0)&=\frac{d}{dz}|_{z=0}\left(-\frac{1}{\Gamma(z)}\int_0^\infty u^{z-1}\Str(Ne^{-\frac{u}{2p}\frD_p^2}) \,du\right)\\
&=\sum_{i\geq 0}\frac{d}{dz}|_{z=0}\left(-\frac{1}{\Gamma(z)}\int_0^\infty u^{z-1}\Str(Na_{i,u}) \,du\right)p^{n-i}\\
&=\sum_{i\geq 0}\beta_ip^{n-i}
\end{split}
\end{equation}
Now,
\begin{equation}\label{thetapexp2}
\begin{split}
\theta_p(z)&=-\frac{1}{\Gamma(z)}\int_0^1 u^{z-1}\Str(Ne^{-\frac{u}{2p}\frD_p^2}) \,du-\frac{1}{\Gamma(z)}\int_1^\infty u^{z-1}\Str(Ne^{-\frac{u}{2p}\frD_p^2}) \,du\\
&=-\sum_{i\geq 0}\frac{p^{n-i}}{\Gamma(z)}\int_0^1 u^{z-1}\Str(Na_{i,u}) \,du -\frac{1}{\Gamma(z)}\int_1^\infty u^{z-1}\Str(Ne^{-\frac{u}{2p}\frD_p^2}) \,du\\
&=-\sum_{i\geq 0}\frac{p^{n-i}}{\Gamma(z)}\sum_{j\geq -n}\int_0^1 u^{z-1}\Str(Na_i^{[j]})u^j \,du -\frac{1}{\Gamma(z)}\int_1^\infty u^{z-1}\Str(Ne^{-\frac{u}{2p}\frD_p^2}) \,du\\
&=-\sum_{i\geq 0}\left(\frac{p^{n-i}}{\Gamma(z)}\sum_{j\geq -n}\frac{\Str(Na_i^{[j]})}{z+j}\right) -\frac{1}{\Gamma(z)}\int_1^\infty u^{z-1}\Str(Ne^{-\frac{u}{2p}\frD_p^2}) \,du\\
\end{split}
\end{equation}
Notice that $\frac{1}{\Gamma(z)}\int_1^\infty u^{z-1}\Str(Ne^{-\frac{u}{2p}\frD_p^2}) \,du$ is a holomorphic function that vanishes at $z=0.$ This implies that
\begin{equation}\label{thetapexp0}
\begin{split}
\theta_p(0)=-\sum_{i\geq 0}\Str(Na_i^{[0]})p^{n-i}
\end{split}
\end{equation}
Now \eqref{thetapexp}, \eqref{thetapexp2} and \eqref{thetapexp0} imply that
\begin{equation}\label{thetamainexp}
\begin{split}
\theta_p'(0)&=\sum_{i\geq 0}\Str(Na_i^{[0]})p^{n-i}\ln p+\sum_{i\geq 0}\beta_ip^{n-i}\\
&=\sum_{i\geq 0}(\alpha_i\ln p+\beta_i)p^{n-i},
\end{split}
\end{equation}
which proves \eqref{introthm}. More explicitly, we have the following identities 
\begin{equation}\label{prize}
\begin{split}
\alpha_0&=\str (Na_0^{[0]})\\
\beta_0&=\frac{d}{dz}|_{z=0}\big(-\frac{1}{\Gamma(z)}\int_0^\infty u^{z-1}\Str(Na_{0,u}) \,du\big)\\
\alpha_1&=\str (Na_1^{[0]})
\end{split}
\end{equation}
Our main result is a calculation of explicit expressions, in terms of integrals of characteristic classes, for $\alpha_0$, $\beta_0$ and, under additional hypotheses, also $\alpha_1$.
\subsubsection{Computation of $\alpha_0$}
In order to compute $\alpha_0$, we need to isolate the terms in \eqref{pexp1} with the highest $p-$power and use them in \eqref{exp2}. This is done with the aid of the modified Getzler rescaling of definition \ref{pgrad}.
\begin{lemma}\label{nicelemma}
As $p\to\infty$ and $u$ near zero, 
\begin{equation}\label{fund}
p^{-n}e^{-\frac{u}{2p}\frD_p^2}(y,y)=\frac{1}{(2\pi)^n}\dfrac{e^{uw_d}\mathrm{det}(\dot{R}^L)}{\mathrm{det}(1-e^{-u\dot{R}^L})}\otimes \mathrm{id}_E(y)+o(p^{-1/2}),
\end{equation}
in the $C^l$-norm for every $l>0$.
\end{lemma}
\begin{proof}
The contributing terms with the highest possible power of $p$ in $\theta_{j,p}(0)$ can be found among the terms in the asymptotic expansion of the heat kernel of $-(\p_{\rmx_i}+p\Gamma^L)^2+pc(R^L)$ which is the result of removing all terms with non-maximal $p$-grading from \eqref{tdirac2}. The rest of the proof follows as in \cite[thm. 4]{alh}. For a different proof of \eqref{fund} in the case $B=0$ see \cite[thm 5.5.8]{mm}.
\end{proof}
Using \eqref{fund}, we can write $p^{-j}\Theta_{j,p}(y,y)=\kappa_j(y,y)+o(p^{-1/2})$ with the $\kappa_j(y,y)$ independent of $p$ such that
$$\dfrac{\strr(Ne^{u\omega_d})\mathrm{det}(\dot{R}^L)}{\mathrm{det}(1-e^{-u\dot{R}^L})}\otimes \mathrm{id}_E(y)\sim \sum_{j=0}^\infty \frac{u^{j-n}}{2^j}\strr(\kappa_j(y,y)).$$
Since these expansions converge uniformly, we get the following expansion for the super trace
\begin{equation}\label{strexp}
\begin{split}
p^{-n}\Str\big(Ne^{-\frac{u}{2p}\frD_p^2}\big)&\sim \mathrm{rk}(E)\int_M\dfrac{\strr(Ne^{u\omega_d})\mathrm{det}(\dot{R}^L/2\pi)}{\mathrm{det}(1-e^{-u\dot{R}^L})}\,dv_y+o(p^{-1/2})\\
\end{split}
\end{equation}
Denote $\alpha_0=\frac{n}{2}\int_M \frac{\omega^n(y)}{n!}\,dv_y.$ Algebraic manipulations \cite[(5.5.37)-(5.5.40)]{mm} give the following expansion where $\alpha_0^{[-1]}$ is a constant
\begin{equation}\label{strexp2}
\begin{split}
\int_M\dfrac{\strr(Ne^{u\omega_d})\mathrm{det}(\dot{R}^L/2\pi)}{\mathrm{det}(1-e^{-u\dot{R}^L})}\,dv_y&=\frac{\alpha_0^{[-1]}}{u}+\frac{n}{2}\int_M \frac{\omega^n(y)}{n!}\,dv_y+\cO(u)\\
\end{split}
\end{equation}
as $u\to0^+.$ We conclude
\begin{lemma}
\begin{equation}\label{alphazeroexp}
\alpha_0=\str (Na_0^{[0]})=\frac{n\, \rke}{2}\int_M \frac{\omega^n(y)}{n!}\,dv_y=\frac{n\,\rke}{2}\int_M 
\mathrm{det}(\dot{R}^L/2\pi)\,dv_y
\end{equation}
\end{lemma}
\subsubsection{Computation of $\beta_0$}
The computation of $\beta_0$ is identical to the one in \cite{bis2} or \cite[Sec. 5.5]{mm} because the leading terms in \eqref{fund} do not involve $\frB$. Instead of just quoting those sources, we present some details of their calculation.


Define
\begin{equation}\label{auxzeta}
\hat{\zeta}(z)=\frac{-1}{\Gamma(z)}\int_0^\infty u^{z-1}\int_M\dfrac{\strr(Ne^{u\omega_d})\mathrm{det}(\dot{R}^L/2\pi)}{\mathrm{det}(1-e^{-u\dot{R}^L})}\,dv_y\,du
\end{equation}
Using \eqref{thetapexp}, \eqref{thetapexp2}, \eqref{prize}, and {\eqref{fund}}  gives
\begin{equation}\label{almost}
p^{-n}\theta_p(0)\ln p+p^{-n}\theta_p'(0)=\hat{\zeta}'(0) \mathrm{rk}(E)+o(p^{-1/2}),
\end{equation}
which implies that, in order to compute $\beta_0$, we only need to compute $\hat{\zeta}'(0). $

At this point we need the identity \cite[(5.5.45)]{mm}
\begin{equation}\label{newid}
\dfrac{\strr(Ne^{u\omega_d})\mathrm{det}(\dot{R}^L/2\pi)}{\mathrm{det}(1-e^{-u\dot{R}^L})}=\mathrm{det}(\dot{R}^L/2\pi)\tr\big((1-\exp(u\dot{R}^L))^{-1}\big)
\end{equation}
Consider an orthonormal frame $\{w_j\}_{j=1}^n$ of $T^{1,0}M$ in which $$\dot{R}^L=\mathrm{diag}(a_1(y),\ldots, a_n(y)).$$
 The positivity asumption on $L$ implies that $a_j(y)>0$ for all $y\in M$ and all $j=1,\ldots,n.$ We get that
 \begin{equation}\label{aux5}
 \tr(1-e^{u\dot{R}^L})^{-1}=\sum_{j=1}^n\dfrac{-e^{-ua_j}}{1-e^{-ua_j}}.
\end{equation}
Finally, plug \eqref{newid} and \eqref{aux5} into \eqref{auxzeta}, perform an obvious change of variables in the resulting integral, and simplify  to get
\begin{equation}\label{punchline}
\hat{\zeta}(z)=\int_M\mathrm{det}(\dot{R}^L/2\pi)\tr((\dot{R}^L)^{-z})\frac{1}{\Gamma(z)}\int_0^\infty u^{z-1}\frac{e^{-u}}{1-e^{-u}}\,du\,dv_y
\end{equation}
The term $\frac{1}{\Gamma(z)}\int_0^\infty u^{z-1}\frac{e^{-u}}{1-e^{-u}}\,du$ equals the usual Riemann $\zeta$-function $\zeta(z)$ and we have $\zeta(0)=-1/2$ and $\zeta'(0)=\frac{-\ln (2\pi)}{2}.$ It follows that
\begin{equation}\label{zetader} 
\beta_0=\rke \hat{\zeta}'(0)=\frac{\rke}{2}\int_M\mathrm{det}(\dot{R}^L/2\pi) \ln \mathrm{det}(\dot{R}^L/2\pi)\,dv_y
\end{equation}
\section{Computation of $\alpha_1$}
The explicit calculation of the $\alpha_1$ term in \eqref{prize} requires considerable work. We start by imposing the following restriction on the positive line bundle $L$.
\begin{assumption}
\begin{equation}\label{assumptionl}
\ome=\frac{\sqrt{-1}}{2\pi}R^L(\cdot,\cdot)=g^{TM}(J\cdot,\cdot)
\end{equation}
\end{assumption}
This implies that $\dot{R}^L=2\pi \mathrm{Id}$ when written in an orthonormal frame $\{w_j\}_{j=1}^n$ of $T^{1,0}M.$
\subsection{Rescaling of Operators}
The main technique used to compute $\alpha_1$ in \eqref{prize} is based on a rescaling trick of Getzler. The function $\kappa$ is defined in \eqref{norcoords}.

\begin{define}
Let $t=p^{-1/2}$ extended to a continuous real parameter taking nonnegative values. Given a smooth section $s\in \Gamma(\cE_p)$ one defines
\begin{equation}\label{tvernew}
\begin{split}
(S_ts)(\rmx)&=s(\rmx/t)\\
\nabla_t&=S_t^{-1}\kappa^{1/2}t\nabla^{\cE_p}\kappa^{-1/2}S_t\\
\caL_2^t&=S_t^{-1}\kappa^{1/2}t^2\frD_p^2\kappa^{-1/2}S_t
\end{split}
\end{equation}
Also, define 
\begin{equation}\label{0nabladef}
\nabla_{0,\p_i}=\p_i+\frac{1}{2}R^L(\frR,\p_i)
\end{equation}
\end{define}
One has that
\begin{equation}\label{nablatcal2}
\begin{split}
\nabla^{\cE_p}_ts(\rmx)&=((d+t\Gamma_{t\rmx}^{S}+tp\Gamma^L_{t\rmx}+t\Gamma_{t\rmx}^{E_{y}})s)(\rmx)-\frac{t}{2\kappa}\p_i(\kappa)(t\rmx)s(\rmx)
\end{split}
\end{equation}
Using  \eqref{0nabladef}, \eqref{gammaexpgen},   \eqref{curvlexp}, and \eqref{gammaLfullexpansion} one gets
\begin{equation}\label{fullnablat2}
\begin{split}
\nabla_{t,\p_i}&=\p_i+t\Gamma_{t\rmx}^{S}+\frac{1}{t}\Gamma_{t\rmx}^L(\p_i)+t\Gamma_{t\rmx}^E(\p_i)-\frac{t}{2\kappa}(\p_i\kappa)(t\rmx)\\
&=(\p_i+\frac{1}{2}R^L_{y}(\frR,\p_i))-\frac{t^2}{6}\big(\pi \langle R^{TM}_{y}(z,\bar{z})\frR,\p_i\rangle_{y} -\langle R^{TM}_{y}(\frR,\p_k)\p_k,\p_i\rangle_{y}\big)\\
&\,\,\,+\frac{t^2}{2}(R^E_{y}(\frR,\p_i)+R_{y}^{S}(\frR,\p_i))+\mathcal{O}(t^3)\\
\end{split}
\end{equation}

Similarly, 
\begin{equation}
\begin{split}
\caL^t_2&=-g^{ij}(t\rmx)\big(\nabla_{t,\p_i}\nabla_{t,\p_j}-t\Gamma_{ij}^l(t\rmx)\nabla_{t,\p_l}\big)\\
&\,\,\,+t^2\big(\frac{r^M}{4}-2|\frB|^2+pc(R^L)+c(R^E+\frac{1}{2}R^{\det})\big)(t\rmx)\\
&=-g^{ij}(t\rmx)\big(\nabla_{t,\p_i}\nabla_{t,\p_j}-t\Gamma_{ij}^l(t\rmx)\nabla_{t,\p_l}\big)+c(R^L)(t\rmx)\\
&\,\,\,+t^2\big(\frac{r^M}{4}-2|\frB|^2+c(R^E+\frac{1}{2}R^{\det})\big)(t\rmx)\\
&=L_2^0+O_1t+O_2t^2+\cO(t^3)
\end{split}
\end{equation}
Expanding each term using \eqref{norcoords} and comparing $t$-degrees one gets

\begin{equation}\label{o2}
\begin{split}
L^0_2&=-(\nabla_{0,\p_i})^2+4\pi N-2\pi n\\
O_1&=0\\
O_2&=\frac{r_{y}^M}{12}-2|\frB|_{y}^2+c(R_{y}^E+\frac{1}{2}R_{y}^{\det})+\frac{1}{3}\langle R^{TM}_{y}(\frR,\p_i)\frR,\p_j\rangle_{y}\nabla_{0,\p_i}\nabla_{0,\p_j}\\
&\,\,+\frac{\pi}{3} \langle R^{TM}_{y}(z,\bar{z})\frR,\p_i\rangle_{y}\nabla_{0,\p_i} -(R^E_{y}(\frR,\p_i)+R_{y}^{S}(\frR,\p_i)))\nabla_{0,\p_i}\\
&\,\,+\frac{1}{3}\langle R^{TM}(\frR,\p_i)\p_i,\p_l\rangle_{y}\nabla_{0,\p_l}\\
\end{split}
\end{equation}

We will need to write the different terms in \eqref{o2} using complex coordinates $\{z_j, \oz_j\}_{j=1,\ldots,n}$ of $\bbR^{2n}=T_yM.$ The radial vector field is then $\frR=z_j\p_{z_j}+\bar{z}_j\p_{\bar{z}_j}=z+\bar{z},$ and $J\frR=iz-i\bar{z}.$ One has that $\sqrt{2}\p_{z_j}=w_j$ for every $j.$ Notice also that 
\begin{equation}\label{complexcoordinates}
\begin{split}
e_{2j-1}&=\p_{z_j}+\p_{\oz_j}\\
e_{2j}&=\sqrt{-1}(\p_{z_j}-\p_{\oz_j})\\
\end{split}
\end{equation}
\begin{define}
The creation and annihilation operators are defined by
\begin{equation}\label{cran}
\begin{split}
b_j&=-2\nabla_{0,\p_{z_j}}=-2\p_{z_j}+\pi \oz_j\\
b_j^+&=2\nabla_{0,\p_{\oz_j}}=2\p_{\oz_j}+\pi z_j
\end{split}
\end{equation}
\end{define}
We use \eqref{cran} and the complex coordinates to rewrite each term in \eqref{o2}. 
\begin{equation}\label{2piterm}
\begin{split}
\frac{\pi}{3} \langle R^{TM}_{y}(z,\bar{z})\frR,\p_i\rangle_{y}\nabla_{0,\p_i}&=\sum_{i\,\,\mathrm{odd}}\frac{\pi}{3} \langle R^{TM}_{y}(z,\bar{z})\frR,\p_i\rangle_{y}\nabla_{0,\p_i}\\
&\,\,\,+\sum_{i\,\, \mathrm{even}}\frac{\pi}{3} \langle R^{TM}_{y}(z,\bar{z})\frR,\p_i\rangle_{x_0}\nabla_{0,\p_i}\\
&=\frac{\pi}{3} \langle R^{TM}_{y}(z,\bar{z})\frR,\p_{2j-1}\rangle_{y}\nabla_{0,\p_{2j-1}}\\
&\,\,\,+\frac{\pi}{3} \langle R^{TM}_{y}(z,\bar{z})\frR,\p_{2j}\rangle_{y}\nabla_{0,\p_{2j}}\\
&=\frac{\pi}{3} \langle R^{TM}_{y}(z,\bar{z})(z+\oz),\p_{z_j}+\p_{\oz_j}\rangle_{y}\nabla_{0,\p_{z_j}+\p_{\oz_j}}\\
&\,\,\,-\frac{\pi}{3} \langle R^{TM}_{y}(z,\bar{z})(z+\oz),\p_{z_j}-\p_{\oz_j}\rangle_{y}\nabla_{0,\p_{z_j}-\p_{\oz_j}}\\
&=\frac{2\pi}{3} \langle R^{TM}_{y}(z,\bar{z})z,\p_{\oz_j}\rangle_{y}\nabla_{0,\p_{z_j}}\\
&\,\,\,+\frac{2\pi}{3} \langle R^{TM}_{y}(z,\bar{z})\oz,\p_{z_j}\rangle_{y}\nabla_{0,\p_{\oz_j}}\\
&=-\frac{\pi}{3} \langle R^{TM}_{y}(z,\bar{z})z,\p_{\oz_j}\rangle_{y}b_j\\
&\,\,\,+\frac{\pi}{3} \langle R^{TM}_{y}(z,\bar{z})\oz,\p_{z_j}\rangle_{y}b_j^+
\end{split}
\end{equation}
For $\bullet=S,E$ one has
\begin{equation}\label{eclterm}
\begin{split}
R_y^\bullet(\frR,\p_i)\nabla_{0,\p_i}&=R_y^\bullet(z+\oz,\p_{2j-1})\nabla_{0,\p_{2j-1}}+R_y^\bullet(z+\oz,\p_{2j})\nabla_{0,\p_{2j}}\\
&=R_y^\bullet(z+\oz,\p_{z_j}+\p_{\oz_j})\nabla_{0,\p_{z_j}+\p_{\oz_j}}\\
&\,\,\,-R_y^\bullet(z+\oz,\p_{z_j}-\p_{\oz_j})\nabla_{0,\p_{z_j}-\p_{\oz_j}}\\
&=2R_y^\bullet(z+\oz,\p_{\oz_j})\nabla_{0,\p_{z_j}}+2R_y^\bullet(z+\oz,\p_{z_j})\nabla_{0,\p_{\oz_j}}\\
&=-R_y^\bullet(z,\p_{\oz_j})b_j+R_y^\bullet(\oz,\p_{z_j})b_j^+
\end{split}
\end{equation}
Another one
\begin{equation}\label{eclterm2}
\begin{split}
\frac{1}{3}\langle R_y^{TM}(\frR,\p_i)\p_i,\p_l\rangle_{y}\nabla_{0,\p_l}&=\frac{1}{3}\langle R_y^{TM}(\frR,\p_i)\p_i,\p_{2j-1}\rangle_{y}\nabla_{0,\p_{2j-1}}\\
&\,\,\,+\frac{1}{3}\langle R_y^{TM}(\frR,\p_i)\p_i,\p_{2j}\rangle_{y}\nabla_{0,\p_{2j}}\\
&=\frac{1}{3}\langle R_y^{TM}(\frR,\p_i)\p_i,\p_{z_j}+\p_{\oz_j}\rangle_{y}\nabla_{0,\p_{z_j}+\p_{\oz_j}}\\
&\,\,\,-\frac{1}{3}\langle R_y^{TM}(\frR,\p_i)\p_i,\p_{z_j}-\p_{\oz_j}\rangle_{y}\nabla_{0,\p_{z_j}-\p_{\oz_j}}\\
&=\frac{2}{3}\langle R_y^{TM}(\oz,\p_{z_k})\p_{\oz_k},\p_{z_j}\rangle_{y}b_j^+\\
&\,\,\,-\frac{2}{3}\langle R_y^{TM}(z,\p_{\oz_k})\p_{z_k},\p_{\oz_j}\rangle_{y}b_j
\end{split}
\end{equation}
Similar calculations give
\begin{equation}\label{ijterm}
\begin{split}
\frac{1}{3}\langle R^{TM}_{y}(\frR,\p_i)\frR,\p_j\rangle_{y}\nabla_{0,\p_i}\nabla_{0,\p_j}&=\frac{1}{3}\langle R^{TM}_{y}(\oz,\p_{z_i})\oz,\p_{z_j}\rangle_{y}b_i^+b_j^+\\
&\,\,\,+\frac{1}{3}\langle R^{TM}_{y}(z,\p_{\oz_i})z,\p_{\oz_j}\rangle_{y}b_ib_j\\
&\,\,\,-\frac{1}{3}\langle R^{TM}_{y}(z,\p_{\oz_i})\oz,\p_{z_j}\rangle_{y}b_ib_j^+\\
&\,\,\,-\frac{1}{3}\langle R^{TM}_{y}(\oz,\p_{z_i})z,\p_{\oz_j}\rangle_{y}b_i^+b_j
\end{split}
\end{equation}
We also need to work with the Clifford product terms in  \eqref{o2}.
\begin{equation}\label{cliffterm}
\begin{split}
c(R_y^\bullet)=\frac{1}{2}R_y^\bullet(\p_l,\p_m)c_lc_m&= 2R_y^\bullet (\p_{z_l},\p_{\oz_m})\big(-\iota_{\bar{w}_l}\varepsilon(\bar{w}^m)+\varepsilon(\bar{w}^m)\iota_{\bar{w}_l}\big)\\
&=\sum_{l\neq m}4R_y^\bullet (\p_{z_l},\p_{\oz_m})\varepsilon(\bar{w}^m)\iota_{\bar{w}_l}\\
&\,\,\,+\sum_l2R_y^\bullet (\p_{z_l},\p_{\oz_l})\big(-\mathrm{Id}+2\varepsilon(\bar{w}^l)\iota_{\bar{w}_l}\big)\\
&=4R_y^\bullet (\p_{z_l},\p_{\oz_m})\varepsilon(\bar{w}^m)\iota_{\bar{w}_l}-2R_y^\bullet (\p_{z_l},\p_{\oz_l})
\end{split}
\end{equation}
From \cite[(1.4.30)]{mm} it follows that
$$R_y^{\det}(\p_{z_l},\p_{\oz_l})=r^M_{y}/4$$
We finally get
\begin{equation}\label{o4b}
\begin{split}
O_2&=\frac{1}{3}\langle R^{TM}_{y}(\oz,\p_{z_i})\oz,\p_{z_j}\rangle_{y}b_i^+b_j^++\frac{1}{3}\langle R^{TM}_{y}(z,\p_{\oz_i})z,\p_{\oz_j}\rangle_{x_0}b_ib_j\\
&-\frac{1}{3}\langle R^{TM}_{y}(z,\p_{\oz_i})\oz,\p_{z_j}\rangle_{y}b_ib_j^+-\frac{1}{3}\langle R^{TM}_{y}(\oz,\p_{z_i})z,\p_{\oz_j}\rangle_{y}b_i^+b_j\\
&-\frac{\pi}{3} \langle R^{TM}_{y}(z,\bar{z})z,\p_{\oz_j}\rangle_{y}b_j+\frac{\pi}{3} \langle R^{TM}_{y}(z,\bar{z})\oz,\p_{z_j}\rangle_{y}b_j^+\\
&+\frac{2}{3}\langle R_y^{TM}(\oz,\p_{z_k})\p_{\oz_k},\p_{z_j}\rangle_{y}b_j^+-\frac{2}{3}\langle R_y^{TM}(z,\p_{\oz_k})\p_{z_k},\p_{\oz_j}\rangle_{y}b_j\\
&+R_y^E(z,\p_{\oz_j})b_j-R_y^E(\oz,\p_{z_j})b_j^+\\
&+R_y^S(z,\p_{\oz_j})b_j-R_y^S(\oz,\p_{z_j})b_j^+\\
&-\frac{r_{y}^M}{6}-2|\frB|_{y}^2\\
&+4R_y^E (\p_{z_l},\p_{\oz_m})\varepsilon(\bar{w}^m)\iota_{\bar{w}_l}-2R_y^E (\p_{z_l},\p_{\oz_l})\\
&+2R_y^{\det} (\p_{z_l},\p_{\oz_m})\varepsilon(\bar{w}^m)\iota_{\bar{w}_l}
\end{split}
\end{equation}
\subsection{Heat Kernel of Rescaled Laplacian}
Recall that $t=p^{-1/2}$. Given a smooth section $s\in \Gamma(M,\cE_p)$, one has that $(e^{-uL^t_2}s)(\rmx)$ equals
\begin{equation}\label{thk1}
\begin{split}
(S^{-1}_t\kappa^{1/2}e^{-\frac{u}{2p}\frD_p^2}&\kappa^{-1/2}S_t s)(\rmx)\\
&=\kappa^{1/2}(t\rmx)\int_M e^{-\frac{u}{p}\frD_p^2}(t\rmx,\rmx')s(\rmx'/t)\kappa^{-1/2}(\rmx')dv(\rmx')\\
&=\kappa^{1/2}(t\rmx)\int_M e^{-\frac{u}{2p}\frD_{2p}^2}(t\rmx,\rmx')s(\rmx'/t)\kappa^{1/2}(\rmx')d\rmx'\\
&=\kappa^{1/2}(t\rmx)\int_M e^{-\frac{u}{2p}\frD_p^2}(t\rmx,tW')s(W')\kappa^{1/2}(tW')t^{2n}dW'\\
&=p^{-n}\kappa^{1/2}(t\rmx)\int_M e^{-\frac{u}{2p}\frD_p^2}(t\rmx,t\rmx')s(\rmx')\kappa^{1/2}(t\rmx')d\rmx'\\
\end{split}
\end{equation}
From \eqref{normalvolume} we have $dv(\rmx')=\kappa(\rmx')d\rmx'$ in normal coordinates $\{\rmx'_j\}.$ Notice that $\kappa(0)=1$. One gets,
\begin{equation}\label{kernelcomp}
\begin{split}
e^{-uL^t_2}(\rmx,\rmx')\kappa(\rmx')&=p^{-n}\kappa^{1/2}(t\rmx)e^{-\frac{u}{2p}\frD_p^2}(t\rmx,t\rmx')\kappa^{1/2}(t\rmx').\\
e^{-\frac{u}{2p}\frD_p^2}(\rmx,\rmx')&=p^ne^{-uL^t_2}({\rmx}/{t},{\rmx'}/{t})\kappa(\rmx'/t)\kappa^{-1/2}(\rmx)\kappa^{-1/2}(\rmx')\\
\end{split}
\end{equation}
A Taylor expansion in $t$ gives
\begin{equation}\label{taylornew}
\begin{split}
e^{-\frac{u}{2p}\frD_p^2}(0,0)&=p^ne^{-uL^t_2}(0,0)\\
e^{-\frac{u}{2p}\frD_p^2}(0,0)&=p^n\sum_{k=0}^\infty \frac{t^k}{k!}(\frac{d^k}{dt^k}e^{-uL^t_2}(0,0))|_{t=0}\\
e^{-\frac{u}{2p}\frD_p^2}(0,0)&=\sum_{k=0}^\infty \frac{p^{n-k/2}}{k!}(\frac{d^k}{dt^k}e^{-uL^t_2}(0,0))|_{t=0}\\
\end{split}
\end{equation}
Comparing the last identity in \eqref{kernelcomp} with \eqref{pexp3} one gets that
\begin{equation}\label{wholepoint}
a_{1,u}=\frac{1}{2}\frac{d^2}{dt^2}e^{-\frac{u}{2}L^t_2}(0,0))|_{t=0}
\end{equation}
\subsubsection{Duhamel Principle}
Let $H_v$ be a smooth family of generalized Laplacians depending on a real parameter $v$, which means operators of the form $H_v=\nabla_v^*\nabla_v+Q_v$ with $Q_v$ a section of an endomorphism bundle, then 
\begin{equation}\label{duhamel}
\p_ve^{-uH_v}=-\int_0^ue^{-(u-u_1)H_v}(\p_vH_v)e^{-u_1H_v}\,du_1
\end{equation}
Furthermore, if $\p_vH_v=0$ then the  chain and product rules imply for the second derivative that
\begin{equation}\label{duhamel2}
\p^2_ve^{-uH_v}=-\int_0^ue^{-(u-u_1)H_v}(\p^2_vH_v)e^{-u_1H_v}\,du_1
\end{equation}
In the case at hand, $H_v=L^t_2$, $\p_t(L^t_2)|_{t=0}=0$, and $\p^2_t(L^t_2)|_{t=0}=2O_2$, so we get
\begin{equation}\label{duhamelt}
\frac{1}{2}\frac{d^2}{dt^2}(e^{-\frac{u}{2}L^t_2}(0,0))|_{t=0}=-\int_0^{u/2}\int_{T_{x_0}M}e^{-(u-u_1)L^0_2}(0,\rmx)O_2 e^{-u_1L^0_2}(\rmx,0)\,du_1
\end{equation}
Here we can identify $\bbR^{2n}\sim T_{x_0}M.$ Therefore, in order to compute $a_{1,u}$, we need to compute $e^{-uL_2^0}(\rmx,0)$ and $O_2 e^{-u_1L^0_2}(\rmx,0).$
\subsubsection{Mehler's Formula}
The heat kernel of $L^0_2(\rmx,0)$ can be written explicitly using the classical Mehler's formula. First define two auxiliary functions
\begin{equation}\label{auxftndef}
\begin{split}
C_{u}&=(1-e^{-2\pi u})^{-n}\\
B_u&=\frac{\pi}{2}\frac{1}{\tanh (\pi u)}=\frac{\pi}{2}\frac{e^{2\pi u}+1}{e^{2\pi u}-1}\\
\end{split}
\end{equation}
The heat kernel is then given by
\begin{lemma}
\begin{equation}\label{findefs}
\begin{split}
e^{-\frac{u}{2}L_{2}^0}(\rmx,0)&=e^{-2\pi u N} C_{u}e^{-B_{u}\|\rmx\|^2}\mathrm{Id}_{\Lambda^{0,*}\otimes E}\\
\end{split}
\end{equation}
\end{lemma}
We will need the following identities
\begin{lemma}
\begin{equation}\label{integrals}
\begin{split}
&\int_0^{u/2}\int_{\bbR^{2n}} e^{-vL_{2}^0}(0,\rmx)e^{-(\frac{u}{2}-v)L_{2}^0}(\rmx,0)\,d\rmx dv=\frac{u}{2}\dfrac{e^{-2u\pi N}}{(1-e^{-2\pi u})^{n}}\\
&\int_0^{u/2}\int_{\bbR^{2n}} e^{-vL_{2}^0}(0,\rmx)|z_j|^2e^{-(\frac{u}{2}-v)L_{2}^0}(\rmx,0)\,d\rmx dv\\
&\,\,\,\,\,\,\,\,=\dfrac{e^{-2u\pi N}}{\pi(1-e^{-2\pi u})^{n+1}}(\frac{u}{2}+\frac{u}{2}e^{-2\pi u}-\frac{1}{2\pi}(1-e^{-2\pi u}))\\
&\int_0^{u/2}\int_{\bbR^{2n}} e^{-vL_{2}^0}(0,\rmx)|z_i|^2|z_j|^2e^{-(\frac{u}{2}-v)L_{2,x}^0}(\rmx,0)\,d\rmx dv\\
&\,\,\,\,\,\,\,\,=\dfrac{2^{\delta_{ij}}e^{-2u\pi N}}{\pi^2(1-e^{-2\pi u})^{n+2}}\big(\frac{u}{2}(1+4e^{-2\pi u}+e^{-4\pi u})-\frac{3}{4\pi}(1-e^{-4\pi u})\big)\\
&b_ie^{-\frac{u}{2}L^0_2}(\rmx,0)=(\pi +2B_u)\bar{z}_ie^{-\frac{u}{2}L^0_2}(\rmx,0)\\
&b_i^+e^{-\frac{u}{2}L^0_2}(\rmx,0)=(\pi-2B_u)z_ie^{-\frac{u}{2}L^0_2}(\rmx,0)\\
\end{split}
\end{equation}
\end{lemma}
\subsection{Computation of $a_{1,u}$}
In this section we use \eqref{wholepoint}, \eqref{duhamelt} and \eqref{integrals} to compute $a_{1,u}.$ When computing \eqref{duhamelt}, we will get several integrals of the form
$$\int_{\bbR^{2n}}P(z_1,\oz_1,\ldots,z_n,\oz_n)e^{-\|z\|^2}dz,$$
where $P$ denotes a polynomial expression. These integrals will vanish unless $P$ is dependent only on the magnitudes $\{|z_j|\}_{j=1}^n.$ In the following expressions, we only record those terms that contribute to \eqref{duhamelt}, and use lower dots at the end to denote the unwritten terms that won't contribute.

The symmetries of the Riemannian curvature tensor can be written in simplified notation as follows
\begin{equation}\label{rsym}
\begin{split}
\langle R^{TM}_{y}(\p_{z_i},\p_{\oz_i})\p_{z_j},\p_{\oz_j}\rangle_{y}=R_{i\oi j\oj}&=R_{i\oj j \oi}=R_{j\oj i\oi}=R_{j\oi i\oj}\\
R_{i\oj k\overline{l}}&=R_{k \overline{l}i\oj}
\end{split}
\end{equation}
Similarly, we sometimes simplify notation by writing
\begin{equation*}
\begin{split}
R_y^E(\p_{z_j},\p_{\oz_j})&=R^E_{j\oj}\\
\langle R^{TM}_{y}(\oz,\p_{z_i})\oz,\p_{z_j}\rangle_{y}&=R_{\oz i\oz j}\,\,\,\mathrm{etc...}
\end{split}
\end{equation*}
We will go term by term in \eqref{o4b} to compute $O_2e^{-\frac{u}{2}L^0_2}$. We start with the first line
\begin{equation}\label{kc0}
\begin{split}
\frac{1}{3}\langle R^{TM}_{y}&(\oz,\p_{z_i})\oz,\p_{z_j}\rangle_{y}b_i^+b_j^+e^{-\frac{u}{2}L^0_2}(\rmx,0)\\
&=\frac{(\pi-2B_u)^2}{3}\langle R^{TM}_{y}(\oz,\p_{z_i})\oz,\p_{z_j}\rangle_{y}z_iz_je^{-\frac{u}{2}L^0_2}(\rmx,0)\ldots\\
&=\frac{(2-\delta_{ij})(\pi-2B_u)^2}{3}\langle R^{TM}_{y}(\p_{\oz_j},\p_{z_i})\p_{\oz_i},\p_{z_j}\rangle_{y}|z_i|^2|z_j|^2e^{-\frac{u}{2}L^0_2}(\rmx,0)\ldots\\
&=\frac{(2-\delta_{ij})(\pi-2B_u)^2}{3}R_{i\oi j\oj}|z_i|^2|z_j|^2e^{-\frac{u}{2}L^0_2}(\rmx,0)\ldots\\
\frac{1}{3}\langle R^{TM}_{y}&(z,\p_{\oz_i})z,\p_{\oz_j}\rangle_{y}b_ib_je^{-\frac{u}{2}L^0_2}(\rmx,0)\\
&=\frac{(\pi+2B_u)^2}{3}\langle R^{TM}_{y}(z,\p_{\oz_i})z,\p_{\oz_j}\rangle_{y}\oz_i\oz_je^{-\frac{u}{2}L^0_2}(\rmx,0)\ldots\\
&=\frac{(2-\delta_{ij})(\pi+2B_u)^2}{3}\langle R^{TM}_{y}(\p_{z_i},\p_{\oz_i})\p_{z_j},\p_{\oz_j}\rangle_{y}|z_i|^2|z_j|^2e^{-\frac{u}{2}L^0_2}(\rmx,0)\ldots\\
&=\frac{(2-\delta_{ij})(\pi+2B_u)^2}{3}R_{i\oi j\oj}|z_i|^2|z_j|^2e^{-\frac{u}{2}L^0_2}(\rmx,0)\ldots\\
\end{split}
\end{equation}
Now the second line
\begin{equation}\label{kc1}
\begin{split}
-\frac{1}{3}\langle R^{TM}_{y}&(z,\p_{\oz_i})\oz,\p_{z_j}\rangle_{y}b_ib_j^+e^{-\frac{u}{2}L^0_2}(\rmx,0)\\
&=\big(-\frac{1}{3}(1-\delta_{ij})R_{z \oi \oz j}b_ib_j^+-\frac{1}{3}R_{z \oi \oz i}b_ib_i^+\big)e^{-\frac{u}{2}L^0_2}(\rmx,0)\ldots\\
&=-\frac{1-\delta_{ij}}{3}R_{z \oi \oz j}(\pi-2B_u)(\pi+2B_u)\oz_iz_je^{-\frac{u}{2}L^0_2}(\rmx,0)\\
&\,\,\,-\frac{1}{3}R_{z\oi \oz i}(\pi-2B_u)(-2+(\pi+2B_u)|z_i|^2)e^{-\frac{u}{2}L^0_2}(\rmx,0)\ldots\\
&=\frac{1-\delta_{ij}}{3}R_{i \oi j \oj}(\pi^2-4B^2_u)|z_i|^2|z_j|^2e^{-\frac{u}{2}L^0_2}(\rmx,0)\\
&\,\,\,+\frac{1}{3}R_{i\oi j\oj}(\pi-2B_u)(-2+(\pi+2B_u)|z_i|^2)|z_j|^2e^{-\frac{u}{2}L^0_2}(\rmx,0)\ldots\\
-\frac{1}{3}\langle R^{TM}_{y}&(\oz,\p_{z_i})z,\p_{\oz_j}\rangle_{y}b_i^+b_je^{-\frac{u}{2}L^0_2}(\rmx,0)\\
&=\big(-\frac{1}{3}(1-\delta_{ij})R_{\oz i z \oj}b^+_ib_j-\frac{1}{3}R_{\oz iz\oi}b^+_ib_i\big)e^{-\frac{u}{2}L^0_2}(\rmx,0)\ldots\\
&=-\frac{1-\delta_{ij}}{3}R_{\oz i z \oj}(\pi+2B_u)(\pi-2B_u)\oz_j z_ie^{-\frac{u}{2}L^0_2}(\rmx,0)\\
&\,\,\,-\frac{1}{3}R_{\oz i z \oi}(\pi+2B_u)(2\delta_{ij}+|z_i|^2(\pi-2B_u))e^{-\frac{u}{2}L^0_2}(\rmx,0)\ldots\\
&=\frac{1-\delta_{ij}}{3}R_{i\oi j \oj}(\pi^2-4B^2_u)|z_i|^2|z_j|^2e^{-\frac{u}{2}L^0_2}(\rmx,0)\\
&\,\,\,+\frac{1}{3}R_{i\oi j\oj}(\pi+2B_u)|z_j|^2(2+|z_i|^2(\pi-2B_u))e^{-\frac{u}{2}L^0_2}(\rmx,0)\ldots\\
\end{split}
\end{equation}
A similar calculation gives the third line. 
\begin{equation}\label{kc2}
\begin{split}
-\frac{\pi}{3} \langle R^{TM}_{y}&(z,\bar{z})z,\p_{\oz_j}\rangle_{y}b_je^{-\frac{u}{2}L^0_2}(\rmx,0)\\
&=-\frac{\pi}{3} R_{z \oz z \oj}\oz_j (\pi+2B_u)e^{-\frac{u}{2}L^0_2}(\rmx,0)\ldots\\
&=-\frac{\pi (\pi+2B_u)}{3} \big(R_{i\oi j\oj}+R_{j \oi i \oj}\big)|z_i|^2|z_j|^2e^{-\frac{u}{2}L^0_2}(\rmx,0)\ldots\\
&=-\frac{(2-\delta_{ij})\pi (\pi+2B_u)}{3} R_{i\oi j\oj}|z_i|^2|z_j|^2e^{-\frac{u}{2}L^0_2}(\rmx,0)\ldots\\
\frac{\pi}{3} \langle R^{TM}_{y}&(z,\bar{z})\oz,\p_{z_j}\rangle_{y}b_j^+e^{-\frac{u}{2}L^0_2}(\rmx,0)\\
&=\frac{\pi(\pi-2B_u)}{3}R_{z \oz \oz j}z_je^{-\frac{u}{2}L^0_2}(\rmx,0)\ldots\\
&=\frac{\pi(\pi-2B_u)}{3}\big(R_{i \oi \oj j}+R_{i \oj \oi j}\big)|z_i|^2|z_j|^2e^{-\frac{u}{2}L^0_2}(\rmx,0)\ldots\\
&=-\frac{(2-\delta_{ij})\pi(\pi-2B_u)}{3}R_{i\oi j\oj}|z_i|^2|z_j|^2e^{-\frac{u}{2}L^0_2}(\rmx,0)\ldots
\end{split}
\end{equation}
The last two terms that involve $R_y^{TM}$ are the following
\begin{equation}\label{kc3}
\begin{split}
\frac{2}{3}\langle R_y^{TM}&(\oz,\p_{z_k})\p_{\oz_k},\p_{z_j}\rangle_{y}b_j^+e^{-\frac{u}{2}L^0_2}(\rmx,0)\\
&=\frac{2(\pi-2B_u)}{3}R_{i\oi j\oj}|z_j|^2e^{-\frac{u}{2}L^0_2}(\rmx,0)\ldots\\
-\frac{2}{3}\langle R_y^{TM}&(z,\p_{\oz_k})\p_{z_k},\p_{\oz_j}\rangle_{y}b_je^{-\frac{u}{2}L^0_2}(\rmx,0)\\
&=-\frac{2(\pi+2B_u)}{3}R_{i\oi j\oj}|z_j|^2e^{-\frac{u}{2}L^0_2}(\rmx,0)\ldots
\end{split}
\end{equation}
The total contribution from the terms involving $R_y^{TM}$ is given by
\begin{equation}\label{partialresult1}
\frac{2-\delta_{ij}}{3}2\pi^2R_{i\oi j\oj}|z_i|^2|z_j|^2e^{-\frac{u}{2}L^0_2}(\rmx,0)\ldots
\end{equation}
For  $\bullet=E,S$ we get
\begin{equation}\label{kc4}
\begin{split}
(R_y^\bullet(z,\p_{\oz_j})b_j-&R_y^\bullet(\oz,\p_{z_j})b_j^+)e^{-\frac{u}{2}L^0_2}(\rmx,0)
\\
&=(R^\bullet_{j\oj}|z_j|^2(\pi+2B_u)-R^\bullet_{\oj j}|z_j|^2(\pi-2B_u))e^{-\frac{u}{2}L^0_2}(\rmx,0)\ldots
\\
&=2\pi R^\bullet_{j\oj}|z_j|^2e^{-\frac{u}{2}L^0_2}(\rmx,0)\ldots,
\end{split}
\end{equation}

Notice that
\begin{equation}\label{scalcomplex}
r^M_{y}=8R_{i\oi j\oj},
\end{equation}
Putting it all together we get
\begin{equation}\label{milestone1}
\begin{split}
O_2e^{-\frac{u}{2}L^0_2}(\rmx,0)&=\frac{2-\delta_{ij}}{3}2\pi^2R_{i\oi j\oj}|z_i|^2|z_j|^2e^{-\frac{u}{2}L^0_2}(\rmx,0)\\
&\,\,\,-\frac{4}{3}R_{i\oi j\oj}e^{-\frac{u}{2}L^0_2}(\rmx,0)\\
&\,\,\,-2R^E_{i\oi}e^{-\frac{u}{2}L^0_2}(\rmx,0)-2|\frB|_{y}^2e^{-\frac{u}{2}L^0_2}(\rmx,0)\\
&\,\,\,+\big(4R_y^E (\p_{z_l},\p_{\oz_m})\varepsilon(\bar{w}^m)\iota_{\bar{w}_l}\\
&\,\,\,+2R_y^{\det} (\p_{z_l},\p_{\oz_m})\varepsilon(\bar{w}^m)\iota_{\bar{w}_l}\big)e^{-\frac{u}{2}L^0_2}(\rmx,0)\\
&\,\,\,+2\pi \big(R^E_{j\oj}+ R^S_{j\oj}\big)|z_j|^2e^{-\frac{u}{2}L^0_2}(\rmx,0)\ldots
\end{split}
\end{equation}
We now use \eqref{wholepoint}, \eqref{duhamelt}, \eqref{integrals}, and \eqref{milestone1} to compute
\begin{equation}\label{progress}
\begin{split}
-&\int_0^{u/2}\int_{\bbR^{2n}}e^{-(u/2-u_1)L^0_2}(0,\rmx)O_2 e^{-u_1L^0_2}(\rmx,0)\,d\rmx\,du_1=\\
&2(R^E_{i\oi}+|\frB|^2_{y})\int_0^{u/2}\int_{\bbR^{2n}}e^{-(u/2-u_1)L^0_2}(0,\rmx) e^{-u_1L^0_2}(\rmx,0)\,d\rmx\,du_1\\
&-2(2R^E_{l\overline{m}}+R^{\det}_{l\overline{m}})\varepsilon(\bar{w}^m)\iota_{\bar{w}_l}\int_0^{u/2}\int_{\bbR^{2n}}e^{-(u/2-u_1)L^0_2}(0,\rmx)e^{-u_1L^0_2}(\rmx,0)\,d\rmx\,du_1\\ 
&+\frac{4}{3}R_{i\oi j\oj}\int_0^{u/2}\int_{\bbR^{2n}}e^{-(u/2-u_1)L^0_2}(0,\rmx)e^{-u_1L^0_2}(\rmx,0)\,d\rmx\,du_1\\
&-2\pi\big(R^E_{j\oj}+ R^S_{j\oj}\big)\int_0^{u/2}\int_{\bbR^{2n}}e^{-(u/2-u_1)L^0_2}(0,\rmx) |z_j|^2e^{-u_1L^0_2}(\rmx,0)\,d\rmx\,du_1\\
&-\frac{2-\delta_{ij}}{3}2\pi^2R_{i\oi j\oj}\int_0^{u/2}\int_{\bbR^{2n}}e^{-(u/2-u_1)L^0_2}(0,\rmx)|z_i|^2|z_j|^2e^{-u_1L^0_2}(\rmx,0)\,d\rmx\,du_1
\end{split}
\end{equation}
Using \eqref{integrals} we get
\begin{equation}\label{progress2}
\begin{split}
-&\int_0^{u/2}\int_{\bbR^{2n}}e^{-(u/2-u_1)L^0_2}(0,\rmx)O_2 e^{-u_1L^0_2}(\rmx,0)\,d\rmx\,du_1=\\
&(R^E_{i\oi}+|\frB|^2_{x_0}-(2R^E_{l\overline{m}}+R^{\det}_{l\overline{m}})\varepsilon(\bar{w}^m)\iota_{\bar{w}_l}+\frac{2}{3}R_{i\oi j\oj})ue^{-2\pi u N}C_{u}\\ 
&-2\pi\big(R^E_{j\oj}+ R^S_{j\oj}\big)\dfrac{e^{-2u\pi N}}{\pi(1-e^{-2\pi u})^{n+1}}(\frac{u}{2}+\frac{u}{2}e^{-2\pi u}-\frac{1}{2\pi}(1-e^{-2\pi u}))\\
&-\frac{2-\delta_{ij}}{3}2\pi^2R_{i\oi j\oj}\dfrac{2^{\delta_{ij}}e^{-2u\pi N}}{\pi^2(1-e^{-2\pi u})^{n+2}}\big(\frac{u}{2}(1+4e^{-2\pi u}+e^{-4\pi u})-\frac{3}{4\pi}(1-e^{-4\pi u})\big)
\end{split}
\end{equation}
Simplify a bit
\begin{equation}\label{progress3}
\begin{split}
-&\int_0^{u/2}\int_{\bbR^{2n}}e^{-(u/2-u_1)L^0_2}(0,\rmx)O_2 e^{-u_1L^0_2}(\rmx,0)\,d\rmx\,du_1=\\
&(R^E_{i\oi}+|\frB|^2_{y}-(2R^E_{l\overline{m}}+R^{\det}_{l\overline{m}})\varepsilon(\bar{w}^m)\iota_{\bar{w}_l}+\frac{2}{3}R_{i\oi j\oj})ue^{-2\pi u N}C_{u}\\ 
&-2\big(R^E_{j\oj}+ R^S_{j\oj}\big)\dfrac{e^{-2u\pi N}}{(1-e^{-2\pi u})^{n+1}}(\frac{u}{2}+\frac{u}{2}e^{-2\pi u}-\frac{1}{2\pi}(1-e^{-2\pi u}))\\
&-\frac{4}{3}R_{i\oi j\oj}\dfrac{e^{-2u\pi N}}{(1-e^{-2\pi u})^{n+2}}\big(\frac{u}{2}(1+4e^{-2\pi u}+e^{-4\pi u})-\frac{3}{4\pi}(1-e^{-4\pi u})\big)
\end{split}
\end{equation}
Notice that $(2-\delta_{ij})2^{\delta_{ij}}=2$ regardless of the values of $i$ and $j$. Finally,
\begin{equation}\label{progress4}
\begin{split}
-&\int_0^{u/2}\int_{\bbR^{2n}}e^{-(u/2-u_1)L^0_2}(0,\rmx)O_2 e^{-u_1L^0_2}(\rmx,0)\,d\rmx\,du_1=\\
&\big[(R^E_{i\oi}+|\frB|^2_{y}-(2R^E_{l\overline{m}}+R^{\det}_{l\overline{m}})\varepsilon(\bar{w}^m)\iota_{\bar{w}_l}+\frac{2}{3}R_{i\oi j\oj})u\\ 
&-2\big(R^E_{j\oj}+ R^S_{j\oj}\big)(1-e^{-2\pi u})^{-1}(\frac{u}{2}+\frac{u}{2}e^{-2\pi u}-\frac{1}{2\pi}(1-e^{-2\pi u}))\\
&-\frac{4}{3}R_{i\oi j\oj}(1-e^{-2\pi u})^{-2}(\frac{u}{2}(1+4e^{-2\pi u}+e^{-4\pi u})-\frac{3}{4\pi}(1-e^{-4\pi u}))\big]\dfrac{e^{-2\pi u N}}{(1-e^{-2\pi u})^{n}}
\end{split}
\end{equation}

\subsubsection{Traces}
What we are really after is $\str (Na_{1,u}).$ To compute this, we will need the following  identities proved in \cite[Lemma 4.6]{fin}. 
\begin{lemma}
\begin{equation}\label{traces}
\begin{split}
\mathrm{str}(Ne^{-u 2\pi N})&=-\rke n e^{-2\pi u}(1-e^{-2\pi u})^{n-1}\\
\mathrm{str}(N (R^E_{l\overline{m}}\varepsilon(\bar{w}^m)\iota_{\bar{w}_l})e^{-u 2\pi N})&=-\mathrm{tr}(R^E_{k\overline{k}})  e^{-2\pi u}(1-ne^{-2\pi u})(1-e^{-2\pi u})^{n-2}\\
\mathrm{str}(N (R^{\det}_{l\overline{m}}\varepsilon(\bar{w}^m)\iota_{\bar{w}_l})e^{-u 2\pi N})&=-\rke R^{\det}_{k\overline{k}}  e^{-2\pi u}(1-ne^{-2\pi u})(1-e^{-2\pi u})^{n-2}
\end{split}
\end{equation}
\end{lemma}
Using \eqref{traces} we get
\begin{equation}\label{str1}
\begin{split}
\mathrm{str}&( Na_{1,u})=-\int_0^{u/2}\int_{\bbR^{2n}}\mathrm{str} Ne^{-(u/2-u_1)L^0_2}(0,\rmx)O_2 e^{-u_1L^0_2}(\rmx,0)\,d\rmx\,du_1=\mathrm{str}N\huge[\\
&(R^E_{i\oi}+|\frB|^2_{y}-(2R^E_{l\overline{m}}+R^{\det}_{l\overline{m}})\varepsilon(\bar{w}^m)\iota_{\bar{w}_l}+\frac{2}{3}R_{i\oi j\oj})u\\ 
&-2\big(R^E_{j\oj}+ R^S_{j\oj}\big)(1-e^{-2\pi u})^{-1}(\frac{u}{2}+\frac{u}{2}e^{-2\pi u}-\frac{1}{2\pi}(1-e^{-2\pi u}))\\
&-\frac{4}{3}R_{i\oi j\oj}(1-e^{-2\pi u})^{-2}(\frac{u}{2}(1+4e^{-2\pi u}+e^{-4\pi u})-\frac{3}{4\pi}(1-e^{-4\pi u}))\huge]\dfrac{ e^{-2\pi u N}}{(1-e^{-2\pi u})^{n}}\\
\end{split}
\end{equation}
which implies
\begin{equation}\label{str2}
\begin{split}
\mathrm{str}&( Na_{1,u})=- n e^{-2\pi u}(1-e^{-2\pi u})^{-1}\huge[\\
&(\mathrm{tr}(R^E_{i\oi})+\rke |\frB|^2_{y}+\frac{2}{3}\rke R_{i\oi j\oj})u\\ 
&-2\mathrm{tr}(R^E_{j\oj})(1-e^{-2\pi u})^{-1}(\frac{u}{2}+\frac{u}{2}e^{-2\pi u}-\frac{1}{2\pi}(1-e^{-2\pi u}))\\
&-\frac{4}{3}\rke R_{i\oi j\oj}(1-e^{-2\pi u})^{-2}(\frac{u}{2}(1+4e^{-2\pi u}+e^{-4\pi u})-\frac{3}{4\pi}(1-e^{-4\pi u}))\huge]\\
&+(\mathrm{tr}(2R^E_{k\overline{k}})+\rke R^{\det}_{k\overline{k}})ue^{-2\pi u}(1-ne^{-2\pi u})(1-e^{-2\pi u})^{-2}\\
&\dfrac{-2}{(1-e^{-2\pi u})^{n}} \mathrm{str} NR^S_{j\oj}e^{-2\pi u N}\big( (1-e^{-2\pi u})^{-1}(\frac{u}{2}+\frac{u}{2}e^{-2\pi u}-\frac{1}{2\pi}(1-e^{-2\pi u}))\big)
\end{split}
\end{equation}
Some cleaning
\begin{equation}\label{str3}
\begin{split}
\mathrm{str}&( Na_{1,u})=
-(\mathrm{tr}(R^E_{i\oi})+\rke |\frB|^2_{y}+\frac{2}{3}\rke R_{i\oi j\oj})\dfrac{nue^{-2\pi u}}{1-e^{-2\pi u}}\\ 
&+2\mathrm{tr}(R^E_{j\oj})\dfrac{ne^{-2\pi u}}{(1-e^{-2\pi u})^{2}}(\frac{u}{2}+\frac{u}{2}e^{-2\pi u}-\frac{1}{2\pi}(1-e^{-2\pi u}))\\
&+\frac{4}{3}\rke R_{i\oi j\oj}\dfrac{ne^{-2\pi u}}{(1-e^{-2\pi u})^{3}}(\frac{u}{2}(1+4e^{-2\pi u}+e^{-4\pi u})-\frac{3}{4\pi}(1-e^{-4\pi u}))\\
&+(\mathrm{tr}(2R^E_{k\overline{k}})+\rke R^{\det}_{k\overline{k}})\dfrac{ue^{-2\pi u}}{(1-e^{-2\pi u})^{2}}(1-ne^{-2\pi u})\\
&-\dfrac{2}{(1-e^{-2\pi u})^{n}} \mathrm{str} NR^S_{j\oj}e^{-2\pi u N}\big( (1-e^{-2\pi u})^{-1}(\frac{u}{2}+\frac{u}{2}e^{-2\pi u}-\frac{1}{2\pi}(1-e^{-2\pi u}))\big)
\end{split}
\end{equation}
It's convenient to isolate the terms with $R^E$ and $|\frB|^2.$
\begin{equation}\label{str4}
\begin{split}
\mathrm{str}&( Na_{1,u})=-\rke |\frB|^2_{y}\dfrac{nue^{-2\pi u}}{1-e^{-2\pi u}}+(\mathrm{tr}R^E_{k\overline{k}})\big(2\frac{u e^{-2\pi u}}{(1-e^{-2\pi u})^2}-\dfrac{ne^{-2\pi u}}{\pi(1-e^{-2\pi u})}\big)\\
&-\frac{2}{3}\rke R_{i\oi j\oj}\dfrac{nue^{-2\pi u}}{1-e^{-2\pi u}}\\
&+\frac{4}{3}\rke R_{i\oi j\oj}\dfrac{ne^{-2\pi u}}{(1-e^{-2\pi u})^{3}}(\frac{u}{2}(1+4e^{-2\pi u}+e^{-4\pi u})-\frac{3}{4\pi}(1-e^{-4\pi u}))\\
&+\rke R^{\det}_{k\overline{k}}\dfrac{ue^{-2\pi u}}{(1-e^{-2\pi u})^{2}}(1-ne^{-2\pi u})\\
&-\dfrac{2}{(1-e^{-2\pi u})^{n}} \mathrm{str} NR^S_{j\oj}e^{-2\pi u N}\big( (1-e^{-2\pi u})^{-1}(\frac{u}{2}+\frac{u}{2}e^{-2\pi u}-\frac{1}{2\pi}(1-e^{-2\pi u}))\big)
\end{split}
\end{equation}

To compute $\alpha_1$ in \eqref{prize}  we need to get the $u=0$ terms in \eqref{str4}. 

Following \cite{fin}, we define the following functions and write down their Laurent expansions
\begin{equation}\label{gdef}
\begin{split}
g_1(u)&=\dfrac{e^{-2\pi u}}{1-e^{-2\pi u}}\\
g_1(u)&=\frac{1}{2\pi u}-\frac{1}{2}+\cO(u)\\
g_2(u)&=\dfrac{e^{-2\pi u}}{(1-e^{-2\pi u})^2}\\
g_2(u)&=g_2^{[-2]}u^{-2}+g_2^{[-1]}u^{-1}-\frac{1}{12}+\cO(u)\\
\tilde{g}_2(u)&=\dfrac{ue^{-2\pi u}}{(1-e^{-2\pi u})^2}\\
\tilde{g}_2(u)&=\tilde{g}_2^{[-1]}u^{-1}+\cO(u)\\
g_3(u)&=\dfrac{ue^{-2\pi u}}{(1-e^{-2\pi u})^3}\\
g_3(u)&=g_3^{[-2]}u^{-2}+g_3^{[-1]}u^{-1}+\cO(u)\\
\end{split}
\end{equation}

Now we collect  lines 2 and 3 in \eqref{str4} and use the following fact
\begin{equation}
R_{i\oi j\oj}=\frac{1}{2}R^{\det}_{j\oj},
\end{equation}
to get 
\begin{equation}\label{str5}
\begin{split}
-\frac{2}{3}&\rke R_{i\oi j\oj}\dfrac{nue^{-2\pi u}}{1-e^{-2\pi u}}\\
&+\frac{4}{3}\rke R_{i\oi j\oj}\dfrac{ne^{-2\pi u}}{(1-e^{-2\pi u})^{3}}(\frac{u}{2}(1+4e^{-2\pi u}+e^{-4\pi u})-\frac{3}{4\pi}(1-e^{-4\pi u}))\\
&=-\frac{2}{3}\rke R_{i\oi j\oj}\dfrac{nue^{-2\pi u}}{1-e^{-2\pi u}}\\
&\,\,\,+\frac{2}{3}\rke R_{i\oi j\oj}\dfrac{ne^{-2\pi u}}{(1-e^{-2\pi u})^{3}}({u}((1-e^{-2\pi u})^2+6e^{-2\pi u})-\frac{3}{2\pi}(1-e^{-4\pi u}))\\
&=4\rke R_{i\oi j\oj}\dfrac{une^{-4\pi u}}{(1-e^{-2\pi u})^{3}}-\dfrac{\rke R_{i\oi j\oj}}{\pi}\dfrac{ne^{-2\pi u}(1+e^{-2\pi u})}{(1-e^{-2\pi u})^2}\\
&=4ne^{-2\pi u}\rke R_{i\oi j\oj}g_3(u)-\dfrac{n\rke R_{i\oi j\oj}}{\pi}g_1(u)-\dfrac{2ne^{-2\pi u}\rke R_{i\oi j\oj}}{\pi}g_2(u)\\
&=n\rke R_{k\overline{k}}^{\det}\big(2e^{-2\pi u}g_3(u)-\dfrac{1}{2\pi}g_1(u)-\dfrac{e^{-2\pi u}}{\pi}g_2(u)\big)
\end{split}
\end{equation}

Going back to the definition of $\nabla^S$ in \eqref{bnabla}, the corresponding curvatures are related by
\begin{equation}\label{curvrel5}
\begin{split}
R^S&=R^{\Lambda^{0,*}}+d_{\nabla^{CL}c(\iota_\bullet\frB)}+ c(\iota_\bullet\frB)\wedge c(\iota_\bullet\frB)\\
&=R^{\Lambda^{0,*}}+\mathfrak{H},
\end{split}
\end{equation}
where $\mathfrak{H}$ is the two-form with values in endomorphisms of $\Lambda^{0,*}T^*M$ defined by the relation. By \cite[(4.25)]{fin}, we have
\begin{equation}\label{curvaux10}
R^{\Lambda^{0,*}}_{j\oj}=R^{\det}_{i\oj}\varepsilon(\bw^j)\iota_{\bw_i}
\end{equation}
Now we can compute a portion of the last supertrace in \eqref{str4}.
\begin{equation}\label{straux}
\begin{split}
&\dfrac{-2}{(1-e^{-2\pi u})^{n}} \mathrm{str}(NR^S_{j\oj}e^{-2\pi u N})\big( (1-e^{-2\pi u})^{-1}(\frac{u}{2}+\frac{u}{2}e^{-2\pi u}-\frac{1}{2\pi}(1-e^{-2\pi u}))\big)\\
&=\dfrac{-2}{(1-e^{-2\pi u})^{n}}\mathrm{str} (NR^{\det}_{i\oj}\varepsilon(\bome^j)\iota_{\bome_i}e^{-2\pi u N})\big((1-e^{-2\pi u})^{-1}(\frac{u}{2}+\frac{u}{2}e^{-2\pi u}-\frac{1}{2\pi}(1-e^{-2\pi u}))\big)\\
&-2\dfrac{\mathrm{str} N\mathfrak{H}_{j\oj}e^{-2\pi u N}}{(1-e^{-2\pi u})^{n}}\big((1-e^{-2\pi u})^{-1}(\frac{u}{2}+\frac{u}{2}e^{-2\pi u}-\frac{1}{2\pi}(1-e^{-2\pi u}))\big)\\
&+2\dfrac{\rke R^{\det}_{j\oj}}{(1-e^{-2\pi u})^{3}}e^{-2\pi u}(1-ne^{-2\pi u})(\frac{u}{2}+\frac{u}{2}e^{-2\pi u}-\frac{1}{2\pi}(1-e^{-2\pi u}))\\
\end{split}
\end{equation}
Going back to \eqref{str4} we have
\begin{equation}\label{strproblem}
\begin{split}
&2\dfrac{\rke R^{\det}_{j\oj}}{(1-e^{-2\pi u})^{3}}e^{-2\pi u}(1-ne^{-2\pi u})(\frac{u}{2}+\frac{u}{2}e^{-2\pi u}-\frac{1}{2\pi}(1-e^{-2\pi u}))\\
&\,\,\,+\dfrac{\rke R^{\det}_{j\oj}}{(1-e^{-2\pi u})^{2}}ue^{-2\pi u}(1-ne^{-2\pi u})\\
&=\dfrac{\rke R^{\det}_{j\oj}(1-ne^{-2\pi u})}{(1-e^{-2\pi u})^2}\big(2ue^{-2\pi u}+\frac{2ue^{-4\pi u}}{(1-e^{-2\pi u})}-\frac{1}{\pi}\big)\\
&=\rke R^{\det}_{j\oj}\big(2g_3(1-ne^{-2\pi u})-\frac{1}{\pi}g_2(1-ne^{-2\pi u})\big)
\end{split}
\end{equation}
Using \eqref{traces}, \eqref{gdef}, \eqref{str5}, \eqref{curvrel5}, and \eqref{straux} one gets
\begin{equation}\label{stralmost}
\begin{split}
\mathrm{str}&( Na_{1,u})=-n\rke |\frB|^2_{y}ug_1(u)+(\mathrm{tr}R^E_{k\overline{k}})\big(\tilde{g}_2(u)-\dfrac{n }{\pi}g_1(u)\big)\\
&\,\,\,+n\, \rke R_{k\overline{k}}^{\det}\big(2e^{-2\pi u}g_3(u)-\dfrac{1}{2\pi}g_1(u)-\dfrac{e^{-2\pi u}}{\pi}g_2(u)\big)\\
&\,\,\,+\rke R^{\det}_{j\oj}\big(2g_3(1-ne^{-2\pi u})-\frac{1}{\pi}g_2(1-ne^{-2\pi u})\big)\\
&\,\,\,-2\dfrac{\mathrm{str} N\mathfrak{H}_{j\oj}e^{-2\pi u N}}{(1-e^{-2\pi u})^{n}}\big((1-e^{-2\pi u})^{-1}(\frac{u}{2}+\frac{u}{2}e^{-2\pi u}-\frac{1}{2\pi}(1-e^{-2\pi u}))\big)\\
\end{split}
\end{equation}
Finally, 
\begin{equation}
\begin{split}
\mathrm{str}&( Na_{1,u})=-\rke |\frB|^2_{y}nug_1(u)+(\mathrm{tr}R^E_{k\overline{k}})\big(\tilde{g}_2(u)-\dfrac{n }{\pi}g_1(u)\big)\\
&\,\,\,+ \rke R_{k\overline{k}}^{\det}\big(2g_3(u)-\dfrac{n}{2\pi}g_1(u)-\frac{1}{\pi}g_2(u)\big)\\
&\,\,\,-2\dfrac{\mathrm{str} N\mathfrak{H}_{j\oj}e^{-2\pi u N}}{(1-e^{-2\pi u})^{n}}\big((1-e^{-2\pi u})^{-1}(\frac{u}{2}+\frac{u}{2}e^{-2\pi u}-\frac{1}{2\pi}(1-e^{-2\pi u}))\big)\\
\end{split}
\end{equation}

Going back to \eqref{prize}, we see that $\alpha_1$ is obtained from the constant terms in the $u-$expansions of each of the functions $ug_1(u), g_2(u),$ and $g_3(u)$ written in \eqref{gdef}. This gives
\begin{equation}\label{alpha1almost}
\begin{split}
&\alpha_1=\int_M \mathrm{str}(Na_1^{[0]})\,dv=\frac{n \rke}{2\pi} \int_M|\frB|^2_{y}dv+\frac{n }{2\pi}\int_M(\mathrm{tr}R^E_{k\overline{k}})\,dv\\
&+ \frac{3n+1}{12\pi}\rke \int_MR_{k\overline{k}}^{\det}\,dv\\
&-2\int_M\big[\dfrac{\mathrm{str} (N\mathfrak{H}_{j\oj}e^{-2\pi u N})}{(1-e^{-2\pi u})^{n}}\big((1-e^{-2\pi u})^{-1}(\frac{u}{2}+\frac{u}{2}e^{-2\pi u}-\frac{1}{2\pi}(1-e^{-2\pi u}))\big)\big]_{u=0}\,dv,\\
\end{split}
\end{equation}
where $\int_M[\cdot]_{u=0}\,dv$ means taking the $u=0$ coefficient. Notice that the volume form can be written in this case as $dv=\frac{\ome^n}{n!}$. Let $\Lambda_\ome$ be the operator defined on forms by $\Lambda_\ome(\eta)\frac{\ome^n}{n!}=\eta\frac{\ome^{n-1}}{(n-1)!}$. We have that $R^{\det}_{j\oj}=\pi \Lambda_\ome(c_1(T^{1,0}M))$ and $R^{E}_{j\oj}=\pi \Lambda_\ome(c_1(E))$, where $c_1(\bullet)$ denotes the first Chern class of the bundle. Using this, we get our final result

\begin{equation}\label{alpha1finalform}
\begin{split}
\alpha_1&=\int_M \mathrm{str}(Na_1^{[0]})\,dv=\frac{n \rke}{2\pi} \int_M|\frB|^2_{y}dv+\frac{n }{2}\int_Mc_1(E)\frac{\ome^{n-1}}{(n-1)!}\\
&\,\,\,+ \frac{(3n+1)\rke}{12} \int_Mc_1(TM)\frac{\ome^{n-1}}{(n-1)!}\\
&\,\,\,-\int_M\big[\dfrac{\mathrm{str} (N\mathfrak{H}_{j\oj}e^{-2\pi u N})}{(1-e^{-2\pi u})^{n+1}}({u}+{u}e^{-2\pi u}-\frac{1}{\pi}(1-e^{-2\pi u}))\big)\big]_{u=0}\,dv\\
\end{split}
\end{equation}
This concludes the proof of \eqref{main1}.

\section{Appendix}
Here we collect additional calculations that are well-known and needed throughout this work.
\subsection{Local Formulas in Normal Coordinates}

Let $\{\rmx_j\}_{j=1}^{2n}$ be normal coordinates around a point $y\in M.$ Remember that $\{\p_{\rmx_j}\}_{j=1}^{2n}=\{\p_j\}_{j=1}^{2n}$ is orthonormal at $y.$ Let $\{\hat{e}_j\}_{j=1}^{2n}$ be an orthonormal local frame obtained by parallel transport of $\{\p_{\rmx_j}\}_{j=1}^{2n}$ along geodesics emanating from $y.$  Let $\frR=\rmx_j\p_j$ denote the radial vector field. We use  $\rmx$ to denote the local coordinates of a point in the normal chart. 

In these coordinates, one has the following convenient Taylor expansions centered at the point $y\in M$, which corresponds to the origin $\rmx=0.$ For an elementary proof see \cite{stern}.
\begin{equation}\label{norcoords}
\begin{split}
\hat{e}_i&=\p_i-\frac{1}{6}\langle R^{TM}_{y}(\frR,\p_i)\frR,\p_j\rangle_{y}\p_j+\mathcal{O}(|\rmx|^3)\\
\p_i&=\hat{e_i}+\frac{1}{6}\langle R^{TM}_{y}(\frR,\p_i)\frR,\p_j\rangle_{y}\hat{e}_j+\mathcal{O}(|\rmx|^3)\\
g_{ij}(\rmx)=\langle \p_i,\p_j\rangle&=\delta_{ij}+\frac{1}{3}\langle R^{TM}_{y}(\frR,\p_i)\frR,\p_j\rangle_{y} +\mathcal{O}(|\rmx|^3)\\
\kappa(\rmx)=|\mathrm{Det}(g_{ij})|^{1/2}(\rmx)&=1+\frac{1}{6}\langle R^{TM}_{y}(\frR,\p_j)\frR,\p_j\rangle_{x_0}+\mathcal{O}(|\rmx|^3)\\
(\p_i\kappa)(\rmx)&=\frac{1}{3}\langle R^{TM}_{y}(\frR,\p_j)\p_i,\p_j\rangle_{y}+\mathcal{O}(|\rmx|^2)\\
\Gamma_{ij}^l(\rmx)&=\frac{1}{3}\big(\langle R^{TM}_{y}(\frR,\p_j)\p_i,\p_l\rangle_{y}+\langle R^{TM}_{y}(\frR,\p_i)\p_j,\p_l\rangle_{y}\big)\\
&\,\,\,+\mathcal{O}(|\rmx|^3)\\
\end{split}
\end{equation}
The volume form in normal coordinates is given by
\begin{equation}\label{normalvolume}
dv(\rmx)=\kappa(\rmx)\,d\rmx
\end{equation}
Let $(F,\nabla^F)$ be any bundle over $M$ with metric connection $\nabla^F$. Using a synchronous frame on it, defined by parallel transport along geodesics emanating from $y$, one gets that
\begin{equation}\label{gammaexpgen0}
\sum_{|\alpha|=r}(\p^\alpha\Gamma^F)_y(\p_l)\frac{\rmx^\alpha}{\alpha!}=\frac{1}{r+1}\sum_{|\alpha|=r-1}(\p^\alpha R^F)_y(\frR,\p_l)\frac{\rmx^\alpha}{\alpha!}
\end{equation}
In particular,
\begin{equation}\label{gammaexpgen}
\Gamma^F(\p_j)(\rmx)=\frac{1}{2}R^F_{y}(\frR,\p_j)+\mathcal{O}(|\rmx|^2),
\end{equation}
\subsection{Curvature of the Line Bundle $L$}
Assume that the line bundle $L$ satisfies $\frac{\sqrt{-1}}{2\pi}R^L(\cdot,\cdot)=g(J\cdot,\cdot)$. This conditions allows us to compute the first few terms of the Taylor expansion of $R^L$ centered at $y$. 

On $T_{y}M$ we can compare the frames $\{\hat{e}_j\}$ and $\{\p_j\}$. In particular, there are functions $\theta_k^l(\rmx)$ such that

\begin{equation}\label{framechange}
\p_k=\theta_k^l(\rmx)\hat{e}_l
\end{equation}
Now
\begin{equation}\label{curvl}
\begin{split}
\frac{\sqrt{-1}}{2\pi}R^L(\p_k,\p_l)&=\langle J\p_k,\p_l\rangle\\
&= \theta_k^i(\rmx)\theta_l^j(\rmx)\langle J\hat{e}_i,\hat{e}_j\rangle_\rmx\\
\end{split}
\end{equation}
Since $M$ is K\"ahler,  $J$ is covariantly constant which implies that   $\langle J\hat{e}_i,\hat{e}_j\rangle_\rmx=\langle J\p_i,\p_j\rangle_{y}$. Now we expand \eqref{curvl} using \eqref{norcoords} and \eqref{framechange}. Notice that
\begin{equation}\label{bob1}
\begin{split}
\delta_{ki}\frac{1}{6}\langle R^{TM}_{y}(\frR,\p_l)\frR,\p_j\rangle_{y}\langle J\p_i,\p_j\rangle_{y}&=\frac{1}{6}\langle R^{TM}_{y}(\frR,\p_l)\frR,J\p_k\rangle_{y}\\
&=\frac{1}{6}\langle R^{TM}_{y}(\frR,J\p_k)\frR,\p_l\rangle_{y}
\end{split}
\end{equation}
Also, 
\begin{equation}\label{bob2}
\begin{split}
\delta_{jl}\frac{1}{6}\langle R^{TM}_{y}(\frR,\p_k)\frR,\p_i\rangle_{y}\langle J\p_i,\p_j\rangle_{y}&=\frac{1}{6}\langle R^{TM}_{y}(\frR,\p_k)\frR,\p_i\rangle_{x_0}\langle J\p_i,\p_l\rangle_{y}\\
&=-\frac{1}{6}\langle R^{TM}_{y}(\frR,\p_k)\frR,J\p_l\rangle_{y}
\end{split}
\end{equation}
From \eqref{bob1} and \eqref{bob2} it follows that
\begin{equation}\label{curvlexp}
\begin{split}
R^L&(\p_k,\p_l)=-2\pi \sqrt{-1}\langle J\p_k,\p_l\rangle_{y}\\
&-\frac{\pi\sqrt{-1}}{3}\big(\langle R^{TM}_{y}(\frR,J\p_k)\frR,\p_l\rangle_{y}-\langle R^{TM}_{y}(\frR,\p_k)\frR,J\p_l\rangle_{y}\big)+\mathcal{O}(|\rmx|^4)
\end{split}
\end{equation}
Comparing with the Taylor expansion centered at $\rmx=0$ one has
\begin{equation}\label{comparisontaylor}
\begin{split}
\sum_{|\alpha|=2}&(\p^\alpha R^L)_{y}(\p_k,\p_l)\frac{\rmx^\alpha}{\alpha!}\\
&=-\frac{\pi\sqrt{-1}}{3}\big(\langle R^{TM}_{y}(\frR,J\p_k)\frR,\p_l\rangle_{y}-\langle R^{TM}_{y}(\frR,\p_k)\frR,J\p_l\rangle_{y}\big),
\end{split}
\end{equation}

Using the complex coordinates $\{z_j, \bar{z}_j\}_{j=1,\ldots,n}$ defined in \eqref{complexcoordinates}, the radial vector field is given by $\frR=z_j\p_{z_j}+\bar{z}_j\p_{\bar{z}_j}=z+\bar{z}.$ Therefore, $J\frR=iz-i\bar{z}.$ This implies, 
\begin{equation}\label{curvlexp2}
\langle R^{TM}_{y}(\frR,J\frR)\frR,\p_l\rangle_{y}=-2\sqrt{-1}\langle R_{y}^{TM}(z,\bar{z})\frR,\p_l\rangle_{y}
\end{equation}

Notice the following fact  
\begin{lemma}
\begin{equation}\label{cancel}
\nabla_{\p_j}\langle R_{y}^{TM}(z,\bar{z})\frR,\p_j\rangle=0
\end{equation}
\end{lemma}
\begin{proof}
See \cite[(4.1.106)]{mm}.
\end{proof}
One has that \eqref{comparisontaylor}, \eqref{curvlexp2}, and \eqref{cancel} imply
\begin{equation}\label{comparisontaylor2}
\begin{split}
\sum_{|\alpha|=2}(\p^\alpha R^L)_{y}(\frR,\p_l)\frac{\rmx^\alpha}{\alpha!}&=-\frac{\pi\sqrt{-1}}{3}\langle R^{TM}_{y}(\frR,J\frR)\frR,\p_l\rangle_{y}\\
&=-\frac{2\pi}{3}\langle R^{TM}_{y}(z,\oz)\frR,\p_l\rangle_{y}\\
\end{split}
\end{equation}
It follows from \eqref{gammaexpgen0} and \eqref{curvlexp} that
\begin{equation}\label{gammaLfullexpansion}
\Gamma_{t\rmx}^L(\p_l)=\frac{1}{2}R^L_y(\frR,\p_l)-\frac{t^2\pi}{6}\langle R_{y}^{TM}(z,\bar{z})\frR,\p_l\rangle_y+\cO(t^3)
\end{equation}


\end{document}